\documentclass[11pt]{article} 
\usepackage[top=20mm,bottom=20mm, left=20mm, right=20mm]{geometry}

\usepackage{authblk}
\usepackage{enumerate}
\usepackage{amsmath, cite, mathrsfs}
\usepackage{graphicx, color, hyperref}
\usepackage{amssymb, amsthm, makecell}
\usepackage{stmaryrd}

 \hypersetup{
     colorlinks=true,       % false: boxed links; true: colored links
     linkcolor=red,          % color of internal links (change box color with linkbordercolor)
     citecolor=blue,        % color of links to bibliography
     filecolor=magenta,      % color of file links
     urlcolor=cyan           % color of external links
 }

\theoremstyle{plain}
\newtheorem{theorem}{Theorem}[section]
\newtheorem{lemma}[theorem]{Lemma}

\theoremstyle{remark}
\newtheorem{remark}{Remark}[section]

\def\inf{\operatorname{inf}}
\def\sd{\operatorname{sd}}
\def\supp{\operatorname{supp}}
\def\singsupp{\operatorname{singsupp}}

\def\div{\operatorname{div}}
\def\Div{\operatorname{Div}}

\def\curl{\operatorname{curl}}
\def\Curl{\operatorname{Curl}}

\def\tr{\operatorname{tr}}
\def\Lin{\operatorname{Lin}}
\def\Sym{\operatorname{Sym}}
\def\Skw{\operatorname{Skw}}

\def\deg{\operatorname{deg}}
\def\dx{\operatorname{dx}}
\def\da{\operatorname{da}}
\def\dl{\operatorname{dl}}
\def\dy{\operatorname{d\boldsymbol{y}}}

\title{Point Singularities in Incompatible Elasticity}
\author{Animesh Pandey and Anurag Gupta\thanks{ag@iitk.ac.in}}
\date{Department of Mechanical Engineering, Indian Institute of Technology Kanpur, 208016, India\\[2ex]%
    \today
}

\begin{document}
\maketitle

\begin{abstract}
The equations of stress equilibrium and strain compatibility/incompatibility are discussed for fields with point singularities in a planar domain. The sufficiency (or insufficiency) of the smooth maps, obtained by restricting the singular fields to the domain away from the singularities, in completely characterizing the equations of equilibrium and compatibility/incompatibility over the entire domain, is established and illustrated with examples. The uniqueness of the solution to the stress problem of incompatible linear elasticity, allowing for singular fields, is proved. The uniqueness fails when the problem is considered solely in terms of the restricted maps. As applications of our framework, a general stress solution, in response to point supported body force and defect fields, is derived and a generalized notion of the force acting on a defect is developed.
\end{abstract}


\section{Introduction}
Physical quantities on a given domain are usually modelled as pointwise maps (functions) which assign them values in an appropriate space for each point in the domain. The stress in an elastic body, for example, is considered as a map with values in the space of symmetric second-order tensors for every point of the domain. The maps are assumed to obey certain regularity conditions depending on the physical problem under consideration. Our interest is in physical fields which have isolated points of singularity such that they are smooth, away from the discrete set of isolated points, but might not be even continuous at the singular points. In fact, the fields may behave in an unbounded manner in the vicinity of a point of singularity. Such fields can appear in the planar problems of incompatible linear elasticity due to the presence of point supported body forces (or force multipoles) and isolated defects (dislocation, dislocation dipole, disclination, interstitial, etc.) \cite{eshelby1956continuum, kroupa1966dislocation, gurtin1973linear, kunin1983elastic, podio2014elasticity}. Clearly, these fields cannot be modelled as pointwise maps over the entire domain since they remain undefined at finite number of points in the domain. This difficulty is circumvented by treating the singular physical fields as distributions. The space of distributions is a generalization of the space of continuous functions which, among other types of singularities, allows for isolated points of singularity \cite{friedlander1998introduction}. Given a distribution with isolated points of singularity, a unique smooth map, away from the points of singularity, can be obtained by restricting the distribution to the set of all non-singular points in the domain. On the other hand, if it exists, the extension of a smooth map (over the non-singular points) to a distribution on the whole domain is not necessarily unique. The precise relationship between the smooth map and the associated distribution, and their derivatives, can be understood using the notion of scaling degree \cite{brunetti2000microlocal}. The scaling degree of a distribution describes its local behaviour at a point under a rescaling of the local domain. We apply these and related ideas to discuss the equations of stress equilibrium and strain compatibility/incompatibility with singular fields. The central contributions of our work are summarized next.

We consider stress ($\boldsymbol{\sigma}$), strain ($\boldsymbol{E}$), body force ($\boldsymbol{B}$), and incompatibility ($N$) fields as distributions with singular support at a given point $O$ in a planar domain $\Omega \subset \mathbb{R}^2$. Let $\boldsymbol{\sigma}_0$, $\boldsymbol{E}_0$, $\boldsymbol{B}_0$, and $N_0$ be the corresponding smooth maps in $\Omega-O$ obtained by restricting the original fields to $\Omega-O$.  We derive a set of local conditions, in terms of $\boldsymbol{\sigma}$ and $\boldsymbol{B}$, equivalent to the distributional form of the stress equilibrium (Lemma~\ref{fullocalstressequib}). In particular, assumptions on $\boldsymbol{\sigma}$ and $\boldsymbol{B}$ are  established under which the local equilibrium conditions are given completely in terms of $\boldsymbol{\sigma}_0$ and $\boldsymbol{B}_0$ (Lemmas~\ref{stressequiv1} and \ref{stressequiv2}). There can be several cases of interest when these assumptions do not hold; e.g., the stress field $\boldsymbol{\sigma}$ developed in response to a force dipole or an extra matter defect does not satisfy the required conditions. We also derive sets of local conditions, in terms of $\boldsymbol{E}$ and $N$, equivalent to distributional forms of strain compatibility and strain incompatibility relations, both within the context of linearized kinematics (Lemmas~\ref{compstrcomplete} and \ref{incompstrcomplete}). Again, assumptions on $\boldsymbol{E}$ and $N$ are established under which the local compatibility and incompatibility conditions are given completely in terms of $\boldsymbol{E}_0$ and ${N}_0$ (Lemmas~\ref{CompatibilityToLoop} and \ref{IncompatibilityCeasroIntegral}). For instance, the strain incompatibility, arising from an isolated defect at $O$, is completely described in terms of $\boldsymbol{E}_0$ and $N_0$ as long as the defect is either a dislocation or a disclination \cite{van2012distributional}. The strain incompatibility due to other kinds of defects (dipoles, extra-matter, concentrated heat source, etc.) cannot be described through $\boldsymbol{E}_0$ and $N_0$ alone.

We prove the uniqueness of solution for the stress boundary value problem of linear incompatible elasticity (Lemma~\ref{LinearElasticityUniquenessLemma}). The problem is governed by stress equilibrium and strain incompatibility conditions, both in their distributional form, in addition to a smooth traction field specified on the boundary of the planar domain. The fields are allowed to be singular at an interior point $O$ of the domain. In particular, we emphasize the assumptions which allow the uniqueness result to hold even when the problem is stated in terms of smooth maps $\boldsymbol{\sigma}_0$, etc. In doing so, we generalize the work by Sternberg and coauthors \cite{sternberg1955concept, Turteltaub1968OnCL, gurtin1973linear}, where the uniqueness results  were established for stresses satisfying certain scaling assumptions near $O$. These assumptions are in fact equivalent to those which allow us to deal with uniqueness in the restricted context of the smooth maps. Our uniqueness result is, however, independent of any a priori assumption on the nature of the stress singularity at $O$. 

We demonstrate the utility of the distributional approach developed in the paper to derive a general stress solution to the incompatibility problem with point supported body force and incompatibility fields, the latter expressible in terms of point supported defect densities and metric anomalies. The general solution can be used to construct the unique solution to any well posed stress boundary value problem with singular fields. 
We also generalize the notion of the force acting on an isolated defect \cite{eshelby1956continuum, kroupa1966dislocation}.  The generalized force is related to changes in the total energy of the elastic domain with respect to infinitesimal translational and non-translational changes in the defect configuration. We derive explicit forms of the generalized force for the cases of an isolated dislocation, an isolated dislocation dipole, and a centre of dilation.  

The paper is organized as follows. In Section~\ref{prelim}, the required mathematical preliminaries from the theory of singular distributions are presented including the concepts of derivatives, singular support, scaling degree and degree of divergence, and singular fields. The equations for stress equilibrium and strain compatibility/incompatibility are discussed in Sections~\ref{stressequibsec} and \ref{straincompincomp}, respectively. The uniqueness of solution to the stress problem in incompatible elasticity is proved in  Section~\ref{uniqueness}. The generalized notion of force acting on a defect is presented in Section~\ref{forceonadefect}. The paper concludes with some prospects on extending the present work in Section~\ref{conc}.

\section{Preliminaries} \label{prelim}
\subsection{Notation} \label{notation}
 Let $\{ \boldsymbol{e}_1, \boldsymbol{e}_2,\dots,\boldsymbol{e}_n\}$ be a fixed orthonormal basis for $\mathbb{R}^n$.  The inner product of any two vectors $\boldsymbol{u},\boldsymbol{v} \in \mathbb{R}^n$ is given by $\langle \boldsymbol{u},\boldsymbol{v} \rangle = u_i v_i$ (summation is implied over repeated indices), where $u_i=\langle \boldsymbol{u},\boldsymbol{e}_i\rangle$. Let $\{\mathbb{R}^n\}^k $ represent the vector space of $k$-linear maps on $\mathbb{R}^n$. The dyadic product $\boldsymbol{u}_1\otimes \boldsymbol{u}_2\otimes \dots\otimes\boldsymbol{u}_k \in \{\mathbb{R}^n\}^k$ of $\boldsymbol{u}_1,\boldsymbol{u}_2,\dotsc, \boldsymbol{u}_k \in \mathbb{R}^n$ is defined by $\boldsymbol{u}_1\otimes \boldsymbol{u}_2\otimes \dots \otimes\boldsymbol{u}_k (\boldsymbol{v}_1,\boldsymbol{v}_2, \dotsc, \boldsymbol{v}_k)=\langle \boldsymbol{u}_1,\boldsymbol{v}_1\rangle \langle \boldsymbol{u}_2,\boldsymbol{v}_2\rangle \dots \langle \boldsymbol{u}_k,\boldsymbol{v}_k\rangle$, for any $\boldsymbol{v}_1,\boldsymbol{v}_2, \dotsc, \boldsymbol{v}_k \in \mathbb{R}^n.$ The inner product on $\{\mathbb{R}^n\}^k$ is given as $\langle \boldsymbol{u}_1\otimes \boldsymbol{u}_2\otimes \dots \otimes\boldsymbol{u}_k, \boldsymbol{v}_1\otimes \boldsymbol{v}_2\otimes \dots \otimes\boldsymbol{v}_k\rangle=\langle \boldsymbol{u}_1,\boldsymbol{v}_1\rangle \langle \boldsymbol{u}_2,\boldsymbol{v}_2\rangle \dots\langle \boldsymbol{u}_k,\boldsymbol{v}_k\rangle$. The elements of $\{\mathbb{R}^n\}^k$ can be interpreted as a linear map from $\{\mathbb{R}^n\}^{k-k_1}$ to $\{\mathbb{R}^n\}^{k_1}$, for $0\leq k_1 \leq k$, such that $\boldsymbol{u}_1\otimes  \dots \otimes\boldsymbol{u}_k (\boldsymbol{v}_1\otimes \dots\otimes\boldsymbol{v}_{k-k_1})= \langle \boldsymbol{u}_{k_1+1}\otimes  \dots \otimes\boldsymbol{u}_{k}, \boldsymbol{v}_1\otimes  \dots \otimes\boldsymbol{v}_{k-k_1}\rangle \boldsymbol{u}_1 \otimes \dots\otimes\boldsymbol{u}_{k_1}$. Let $\Omega \subset \mathbb{R}^n$ be a  bounded, open, connected set such that the origin $O$ is contained in $\Omega$. Let $B_r \subset \Omega$ represent an open ball of radius $r$, centred at point $O$, with boundary $\partial B_r$.
 For any two sets $A$ and $B$, we use $A \subset B$ and $A-B$ to represent that $A$ is a subset of $B$ and the difference of the sets $A$ and $B$, respectively. In particular, given $O\in \Omega$, $\Omega - O$ denotes the difference of the set $\Omega$ and the singleton set $\{O\}$. 
If an arbitrary point $\boldsymbol{x}\in \mathbb{R}^n$ is expressed in terms of its components as $\boldsymbol{x}=x_i\boldsymbol{e}_i$ then, for any differentiable function on $\mathbb{R}$, say $f$, ${\partial f}/{\partial x_i}$ represents the partial derivative of $f$ with respect to the $i$-th component. The gradient of $f$ is a continuous map $\nabla f :\Omega \to \mathbb{R}^n$ given by $\nabla f=({\partial f}/{\partial x_i}) \boldsymbol{e}_i.$

For $\Omega \subset \mathbb{R}^2$, $\Lin$ denotes the space of linear maps from $\mathbb{R}^2$ to itself with inner product as defined for $\{\mathbb{R}^2\}^2$. The identity map in $\Lin$ is denoted as $\boldsymbol{I}$. The trace of $\boldsymbol{A} \in \Lin$ is given by $\tr \boldsymbol{A} =\langle \boldsymbol{A},\boldsymbol{I} \rangle.$ We use $\Sym$ and $\Skw$ to denote the vector space of symmetric and skew symmetric linear maps from $\mathbb{R}^2$ to itself, respectively. For any $\boldsymbol{v}\in \mathbb{R}^2,$ let $\boldsymbol{e}_3\times \boldsymbol{v}=-\boldsymbol{v}\times \boldsymbol{e}_3=v_2 \boldsymbol{e}_1-v_1 \boldsymbol{e}_2.$ Given a vector valued differentiable map, $\boldsymbol{f}:\Omega \to \mathbb{R}^2$, we introduce $\div \boldsymbol{f}:\Omega \to \mathbb{R}$, given by $\div \boldsymbol{f}=({\partial f_1}/{\partial x_1})+({\partial f_2}/{\partial x_2})$, and $\curl \boldsymbol{f}:\Omega \to \mathbb{R}$, given by $\curl \boldsymbol{f}=({\partial f_2}/{\partial x_1})-({\partial f_1}/{\partial x_2})$. Let $\mathbb{A}\in (\mathbb{R}^2)^4$ be a linear map $\mathbb{A}: \Lin \to \Lin$ defined by $\mathbb{A} (\boldsymbol{v}\otimes \boldsymbol{w})= (\boldsymbol{e}_3 \times\boldsymbol{v})\otimes (\boldsymbol{e}_3 \times\boldsymbol{w})$ for any $\boldsymbol{v}, \boldsymbol{w} \in \mathbb{R}^2$. The bilinear map $\mathbb{A}(\boldsymbol{V},\boldsymbol{W})=\left\langle \mathbb{A} \boldsymbol{V}, \boldsymbol{W} \right\rangle$, for any $\boldsymbol{V},  \boldsymbol{W} \in \mathbb{R}^2\otimes\mathbb{R}^2$, is symmetric, i.e., $\mathbb{A}(\boldsymbol{V},\boldsymbol{W})=\mathbb{A}(\boldsymbol{W},\boldsymbol{V})$. We use $C^{\infty}(\Omega)$, $C^{\infty}(\Omega,\mathbb{R}^2)$, and $C^{\infty}(\Omega,\Lin)$ to represent the space of smooth scalar valued, vector valued, and tensor valued maps, respectively. Given $(x_1,x_2)\in \mathbb{R}^2$, let $r=({x_1}^2 + {x_2}^2)^{1/2}$ and $\theta=\tan^{-1}({x_2}/{x_1})$ be the polar coordinates, and $\{ \boldsymbol{e}_r,\boldsymbol{e}_\theta\}$ be the orthonormal polar basis in $\mathbb{R}^2$.

\subsection{Distributions}
For an open set $\Omega \subset \mathbb{R}^n$ let $\mathcal{D}(\Omega)$ represent the space of compactly supported smooth functions from $\Omega$ to $\mathbb{R}$. 
The space of distributions $\mathcal{D}'(\Omega)$ is the space of all linear continuous functions on $\mathcal{D}(\Omega)$ \cite{friedlander1998introduction}. The space of $\mathbb{R}^2$ valued distributions, represented by $\mathcal{D}'(\Omega,\mathbb{R}^2)$, is such that, for any compactly supported $\mathbb{R}^2$ valued smooth function $\boldsymbol{\phi} \in \mathcal{D}(\Omega,\mathbb{R}^2)$ and $\boldsymbol{T} \in \mathcal{D}'(\Omega,\mathbb{R}^2)$, $\boldsymbol{T}(\boldsymbol{\phi})=T_i (\phi_i)$ (summation is implied over repeated indices). The space of $\Lin$ valued distributions, represented by $\mathcal{D}'(\Omega,\Lin)$, is the dual of the space of  compactly supported $\Lin$ valued smooth functions such that, for $\boldsymbol{\phi} \in \mathcal{D}(\Omega,\Lin)$ and $\boldsymbol{T} \in \mathcal{D}'(\Omega,\Lin)$, $\boldsymbol{T}(\boldsymbol{\phi})=T_{ij} (\phi_{ij})$. 

A function $f:\Omega \to \mathbb{R}$ is said to be locally integrable if it is integrable in any compact subset of $\Omega$. Given any locally integrable function $f: \Omega \to \mathbb{R}$, we can associate with it a distribution $T_f \in \mathcal{D}'(\Omega)$ given by $T_f(\psi)=\int_{\Omega} f\psi \dx$ for all $\psi \in \mathcal{D}(\Omega)$, where $\dx$ represents the Lebesgue measure in $\mathbb{R}^n$. Given a distribution $F \in \mathcal{D}'(\Omega)$, if there exists a locally integrable function $f:\Omega \to \mathbb{R}$ such that $F(\psi)=\int_{\Omega} f\psi \dx$ for all $\psi \in \mathcal{D}(\Omega)$ then $F$ is said to be locally integrable. We call a distribution $T\in \mathcal{D}'(\Omega)$ continuous (or differentiable or smooth) if there exists a continuous (or differentiable or smooth) function $f:\Omega \to \mathbb{R}$ such that $T(\psi)=\int_{\Omega} f\psi \dx$ for all $\psi \in \mathcal{D}(\Omega)$. We say that a distribution $T\in \mathcal{D}'(\Omega)$ is the zero distribution if $T(\psi)=0$ for all $\psi \in \mathcal{D}(\Omega)$.  Given a point $O \in \Omega$, we use $\delta_O \in \mathcal{D}'(\Omega)$ to represent the Dirac measure at $O$, i.e., $\delta_O (\psi)=\psi(O)$ for all $\psi \in \mathcal{D}(\Omega)$.
A sequence of distributions $T_j \in \mathcal{D}'(\Omega)$ is said to converge to $T_0 \in \mathcal{D}'(\Omega)$, as $j\to \infty$, in the sense of distributions if $T_j (\phi) \to T_0 (\phi)$
for all $\phi \in \mathcal{D}(\Omega)$.

\subsection{Derivatives of distributions}

The partial derivative of a distribution $T\in\mathcal{D}'(\Omega)$, $\partial_i T\in\mathcal{D}'(\Omega)$, is defined as
$\partial_i T (\psi)= -T ( {\partial \psi}/{\partial x_i} )$ for all $\psi \in \mathcal{D}(\Omega)$.
If $\alpha \in \mathbb{N}^n$ denotes  an $n$ dimensional multi-index, where $\mathbb{N}$ is the set of non-negative integers, then for $\alpha=(\alpha_1,\alpha_2,\dots,\alpha_n),$ $|\alpha|=\alpha_1 +\alpha_2 +\dots+\alpha_n$, and $T \in \mathcal{D}'(\Omega)$, $\partial^{\alpha} T=\partial^{\alpha_1}_1\partial^{\alpha_2}_2\dots\partial^{\alpha_n}_n T$.
For any sequence such that $T_j\to T_0$, as $j\to \infty$, $\partial_i T_j \to \partial_i T_0 $ and $\partial^\alpha T_j \to \partial^\alpha T_0$, for any multi-index $\alpha$.
 
For $\Omega \subset \mathbb{R}^2$, given any scalar field $T\in \mathcal{D}'(\Omega)$, the gradient of $T$, $\nabla T \in \mathcal{D}'(\Omega,\mathbb{R}^2)$, is given by $\nabla T=(\partial_i T) \boldsymbol{e}_i$. For any vector field $\boldsymbol{v} \in \mathcal{D}'(\Omega,\mathbb{R}^2)$ and any $\psi \in \mathcal{D}(\Omega)$, the divergence of $\boldsymbol{v}$, $\Div \boldsymbol{v}\in \mathcal{D}' (\Omega)$, is given by $\Div \boldsymbol{v} (\psi)=-\boldsymbol{v}(\nabla \psi)$ and the curl of $\boldsymbol{v}$, $\Curl \boldsymbol{v}\in \mathcal{D}' (\Omega)$, is given by $\Curl \boldsymbol{v}(\psi)=-\boldsymbol{v}(\boldsymbol{e}_3\times \nabla \psi)$. For any tensor field $\boldsymbol{V} \in \mathcal{D}'(\Omega,\Lin)$ and any $\boldsymbol{\psi} \in \mathcal{D}(\Omega,\mathbb{R}^2)$, the divergence of $\boldsymbol{V}$, $\Div \boldsymbol{V}\in \mathcal{D}' (\Omega, \mathbb{R}^2)$, is given by $\Div \boldsymbol{V} (\boldsymbol{\psi})=-\boldsymbol{V}(\nabla \boldsymbol{\psi})$ and the curl of $\boldsymbol{V}$, $\Curl \boldsymbol{V}\in \mathcal{D}' (\Omega,\mathbb{R}^2)$, is given by $\langle \Curl \boldsymbol{V},\boldsymbol{a}\rangle=\Curl  (\boldsymbol{V}^T \boldsymbol{a})$, for any fixed $\boldsymbol{a}\in \mathbb{R}^2$.

\subsection{Restriction of distributions}

For an open set $\Omega \subset \mathbb{R}^n$ and  $\omega\subset \Omega$, we can extend $\psi \in \mathcal{D}(\omega)$ to $\overline{\psi} \in \mathcal{D}(\Omega)$ such that, for any $\boldsymbol{x}\in \Omega$,
\begin{equation}
\overline{\psi}(\boldsymbol{x})=\begin{cases}
    \psi(\boldsymbol{x})& \text{if } \boldsymbol{x}\in \omega,\\
    0              & \text{otherwise.}
\end{cases}
\end{equation} 
Given a distribution $T\in \mathcal{D}'(\Omega)$ and an open subset $\omega \subset \Omega$, the restriction of $T$ to $\omega$, $T|_{\omega}\in \mathcal{D}'(\omega)$, is given by $T|_{\omega} (\psi)=T(\overline{\psi})$, for all $\psi \in \mathcal{D}(\omega)$. Restricting distributions to arbitrary open sets containing a common point allows us to localize them at the given point.  We say that a property is local at a point $\boldsymbol{x}\in \Omega$ if it is equivalent for any two arbitrary distributions $T_1, T_2 \in \mathcal{D}'(\Omega)$ for which there exists an open set $\omega \subset \Omega$ where $\boldsymbol{x}\in \omega$ and $T_1|_\omega=T_2|_\omega$.
We note that $\partial^{\alpha} \overline{\psi}=\overline{\partial^{\alpha} {\psi}}$ for any $\psi \in \mathcal{D}(\omega)$ and any multi-index $\alpha \in \mathbb{N}^n$. Consequently, the restriction of the derivative of a distribution $T\in\mathcal{D}'(\Omega)$ to an open subset $\omega\subset \Omega$ is same as the derivative of the restriction of the distribution $T|_\omega \in \mathcal{D}'(\omega)$, i.e., $(\partial^{\alpha} T)|_{\omega}= \partial^{\alpha} (T|_\omega)$. In particular, for $\Omega \subset \mathbb{R}^2$, given any scalar field $T\in \mathcal{D}'(\Omega)$, $\nabla(T|_\omega)=(\nabla T )|_\omega$ and, given any vector field $\boldsymbol{v} \in \mathcal{D}'(\Omega,\mathbb{R}^2)$, $\Div (\boldsymbol{v}|_\omega)=(\Div \boldsymbol{v})|_\omega$ and $\Curl (\boldsymbol{v}|_\omega)=(\Curl \boldsymbol{v})|_\omega$.

\subsection{Support and singular support}

The support of a function $f:\Omega\to \mathbb{R}$ is defined as the closure of the set of all points where the function $f\neq 0$. We say that a distribution $T$ is non-zero at a point $\boldsymbol{x}\in \Omega$ if there exists no open set $\omega \subset \Omega$ such that $\boldsymbol{x} \in \omega$ and $T|_{\omega}=0$. The support of a distribution $T\in \mathcal{D}'(\Omega)$, $\supp(T) \subset \Omega$, is defined as the closure of the set of all points in $\Omega$ where $T$ is non-zero. It can be equivalently defined as the smallest closed set $\omega\subset \Omega$ such that $T|_{\Omega-\omega}=0$. Given an integrable function $f$, the support of $T_f \in \mathcal{D}'(\Omega)$ is same as the support of  $f$. We note that $\supp (\partial^\alpha T) \subset \supp (T)$ for any distribution $T\in \mathcal{D}'(\Omega)$ and any multi-index $\alpha$. We say that a distribution is smooth at $\boldsymbol{x}\in \Omega$ if there exists an open set $\omega \subset \Omega$ such that $\boldsymbol{x} \in \omega$ and $T|_\omega$ is smooth. The smoothness of a distribution at any point $\boldsymbol{x} \in \Omega$ is a local property. The singular support of a distribution $T$, $\singsupp(T) \subset \supp(T)$, is defined as the closure of the set of all points in $\Omega$ where $T$ is not smooth. It is the smallest closed set $\omega\subset \Omega$ such that $T|_{\Omega-\omega}$ is smooth, i.e., there exists a smooth function, $f:\Omega-\singsupp(T)\to \mathbb{R}$, such that $T|_{\Omega-\singsupp(T)}(\psi)=\int_{\Omega-\singsupp(T)} f \psi \dx$ for all $\psi \in \mathcal{D}(\Omega-\singsupp(T))$. We note that $\singsupp (\partial^\alpha T) \subset \singsupp (T)$ for any distribution $T\in \mathcal{D}'(\Omega)$ and any multi-index $\alpha$.

Let $\mathcal{P}:\mathcal{D}'(\Omega)\to \mathcal{D}'(\Omega)$ be a linear differential operator of order $m$ with constant coefficients, i.e., for any $T \in \mathcal{D}'(\Omega)$, $\mathcal{P}(T)=\sum_{k=0}^{m} \left\langle \boldsymbol{P}^k, \nabla^k T \right\rangle$,
where $\boldsymbol{P}^k\in (\mathbb{R}^n)^k$ such that $\boldsymbol{P}^m$ is non-zero and $\nabla^k T = \nabla\nabla \dots (k~\text{times}) \dots \nabla T$. The linear differential operator $\mathcal{P}$ is elliptic if 
\begin{equation}
\left\langle \boldsymbol{P}^m, \boldsymbol{x}\otimes \boldsymbol{x} \cdots (m~\text{times}) \cdots \otimes \boldsymbol{x} \right\rangle \neq 0
\end{equation}
 for any non-zero $\boldsymbol{x}\in \mathbb{R}^n$. For any arbitrary linear differential operator, $\singsupp(\mathcal{P}(T))\subset \singsupp(T)$. Hence, if a distribution $T$ is smooth in an open set $\omega$ then the distribution $\mathcal{P}(T)$ will also be smooth in $\omega$. The converse assertion, i.e., if $\mathcal{P}(T)$ is smooth in $\omega$ then $T$ is smooth in $\omega$, holds true if the differential operator is elliptic, as established through the following lemma.
\begin{lemma}
\label{EllipticRegularityLemma}
\textup{\cite[Theorem~8.6.1]{friedlander1998introduction}} Given a elliptic linear operator with constant coefficients $\mathcal{P}$ and an open set $\Omega\subset \mathbb{R}^n$, $\singsupp(\mathcal{P} (T))= \singsupp(T)$
for any $T\in \mathcal{D}'(\Omega)$.
\end{lemma}

The product of two continuous functions is defined pointwise as a continuous function. The product of two arbitrary distributions is not well defined. Given any two smooth functions $f \in C^\infty (\Omega)$ and $\phi \in \mathcal{D}(\Omega)$, the product is a smooth function $f\phi \in \mathcal{D}(\Omega)$. Given a smooth function $f\in C^\infty (\Omega)$ and a distribution $T\in \mathcal{D}'(\Omega)$, the product is a distribution $fT \in \mathcal{D}'(\Omega)$ such that $fT(\phi)=T(f\phi)$ for all $\phi\in \mathcal{D}(\Omega).$ If $T$ is continuous, i.e., there exists a continuous function $f_1$ such that $T(\phi)=\int_\Omega f_1 \phi \dx$, then $fT(\phi)=\int_\Omega f f_1 \phi \dx$. The product of a smooth function with a distribution therefore generalizes the notion of the product of continuous functions. For any sequence $T_j\in \mathcal{D}'(\Omega)$ and $T_0 \in \mathcal{D}'(\Omega)$ such that $T_j \to T_0$, as $j \to \infty$, the sequence $fT_j \in \mathcal{D}'(\Omega)$ converges to $fT_0 \in \mathcal{D}'(\Omega)$. The derivative of the product of a distribution $T$ with smooth function $f$ follows the Leibniz rule,
$\partial_i (fT)= (\partial_i f)T + f(\partial_i T)$. Finally, we note that $\supp (fT) \subset \supp(T)$ and $\singsupp (fT)\subset \singsupp(T)$.

\subsection{Scaling degree and degree of divergence}
\label{sddeg}
For any $\phi \in \mathcal{D}(\mathbb{R}^n)$ and $\lambda \in (0,\infty)$ let $\phi_\lambda \in \mathcal{D}(\mathbb{R}^n)$ be such that, for $\boldsymbol{x}\in \mathbb{R}^n$,
\begin{equation}
\label{RescaledTestFunction}
\phi_{\lambda}(\boldsymbol{x})=\frac{1}{\lambda^n} \phi\left( \frac{\boldsymbol x}{\lambda}\right).
\end{equation}
Given $T \in \mathcal{D}'(\mathbb{R}^n)$, the rescaled distribution $T_{\lambda}\in \mathcal{D}'(\mathbb{R}^n)$  is given by $T_{\lambda}(\phi)=T(\phi_{\lambda})$.
 For a locally integrable distribution $T \in \mathcal{D}'(\mathbb{R}^n),$ such that $T(\phi)=\int f\phi \dx$, $T_\lambda(\phi)=\int f_\lambda \phi \dx$, where $f_\lambda ({\boldsymbol x})=f(\lambda \boldsymbol{x})$. The scaling degree of the distribution $T$, $\sd (T)$, with respect to origin $O \in \mathbb{R}^n$ is defined as \cite{brunetti2000microlocal}
\begin{equation}
\sd (T)=\inf \{  k \in \mathbb{R} | \lim_{\lambda \to 0} \lambda ^k T_{\lambda}=0\}.
\end{equation} 
The degree of divergence of $T \in \mathcal{D}'(\mathbb{R}^n)$, with respect to $O$, is defined as the difference of the scaling degree of $T$ and the dimension of the space,
  \begin{equation}
 \deg (T)=\sd(T)-n. \label{degdiv}
  \end{equation} 
 The following lemma establishes the local nature of scaling degree or degree of divergence. 
\begin{lemma}
For $T_1 \in \mathcal{D}'(\mathbb{R}^n)$ and $T_2 \in \mathcal{D}'(\mathbb{R}^n)$, such that ${T_1}|_{\omega}={T_2}|_{\omega}$ for any open set $\omega \subset \mathbb{R}^n$ with $O\in \omega$, $\deg(T_1)=\deg(T_2)$.
\begin{proof}
Given any $\phi\in \mathcal{D}(\mathbb{R}^n),$ there exists $\lambda_0>0$ such that $\supp(\phi_\lambda) \subset \omega$ for all $\lambda<\lambda_0.$ Hence $\lim_{\lambda \to 0} \lambda ^k( {T_1}_{\lambda} - {T_2}_{\lambda}) = 0$,
which implies $\deg(T_1)=\deg(T_2).$
\end{proof}
\end{lemma}
We can use the local nature of the degree of distributions in $\mathbb{R}^n$ to extend the notions of scaling degree and degree of divergence to distributions on arbitrary open sets containing $O$. For any $\phi \in \mathcal{D}({B}_r)$ let $\phi_\lambda \in \mathcal{D}({B}_r)$, with $0<\lambda <1$, be given by \eqref{RescaledTestFunction}. Given $T\in \mathcal{D}'(\Omega)$ and ${B}_r \subset \Omega$ we define $({T|_{{B}_r}})_\lambda \in \mathcal{D}'({B}_r)$ as
\begin{equation}
({T|_{{B}_r}})_\lambda (\phi)=({T|_{{B}_r}})(\phi_\lambda)
\end{equation}  
for all $\phi \in \mathcal{D}({B}_r).$  The scaling degree of the distribution $T$, with respect to $O$, is then defined as
\begin{equation}
\sd (T)=\inf \{  k \in \mathbb{R} | \lim_{\lambda \to 0} \lambda ^k ({T|_{{B}_r}})_\lambda=0\}.
\end{equation} 
The degree of divergence of $T$ can be evaluated using \eqref{degdiv}.
Further, both of these can be defined for distributions $T\in \mathcal{D}'(\Omega-O)$ by considering the restriction $T|_{{B}_r - O}$. The scaling degree and the degree of divergence of a distribution will always be mentioned with respect to $O$ unless stated otherwise. The following results will be used throughout \cite{brunetti2000microlocal}:
\begin{enumerate}[i.]
\item  For $T\in \mathcal{D}'(\Omega),$ $\sd(T)\geq \sd(T|_{\Omega-O})$ and $\deg(T)\geq \deg(T|_{\Omega-O})$.
\item For $T_1 \in \mathcal{D}'(\Omega)$, such that $T_1(\phi)=\int_{\Omega} f\phi \dx$ and $f$ is smooth, $\sd(T_1)\leq 0$. Given a smooth function $f$ and a distribution $T\in \mathcal{D}'(\Omega),$ $\sd(fT) \leq \sd(T)$ and $\deg(fT)\leq \deg(T)$.
\item For $T\in \mathcal{D}'(\mathbb{R}^n)$, such that $T_\lambda=\lambda^k T$, $\sd(T)=-k$ and $\deg(T)=-k-n$. For $T\in \mathcal{D}'(\Omega-O)$, such that $T(\phi)=\int_{\Omega-O} f\phi \dx$ and $f(\boldsymbol{x})\leq c_0 |\boldsymbol{x}|^k$, $\sd(T) \leq -k$ and $\deg(T) \leq -k-n$.
\item For $T=\delta_O$, $T_\lambda=\lambda^{-n} \delta_O$, which implies $\sd(\delta_O)=n$ and $\deg(T)=0.$ For $T=\partial^\alpha \delta_O$, where $\alpha$ is a multi-index, $T_\lambda=\lambda^{-n-|\alpha|} \partial^\alpha \delta_O,$ which implies $\sd(T)=|\alpha|+n$ and $\deg(T)=|\alpha|$.
\item For any distribution $T$, $\sd (\partial^\alpha T) \leq \sd(T) + |\alpha|$ and $\deg (\partial^{\alpha} T) \leq \deg(T) + |\alpha|$.
\end{enumerate}

\subsection{Extension of distributions}
Given an open subset $\omega \subset \Omega$ and a distribution $T_0\in \mathcal{D}'(\omega)$, we say that ${T}\in \mathcal{D}'(\Omega)$ is an extension of $T_0$ if ${T}|_{\omega}={T_0}$. Such an extension is in general not unique. The existence and uniqueness of an extension for distributions $T_0 \in \mathcal{D}'(\mathbb{R}^n-O)$ with a given degree of divergence is due to Brunetti and Fredenhagen \cite{brunetti2000microlocal}. The analogous result for distributions $T_0 \in \mathcal{D}'(\Omega-O)$ is stated in the following lemma. The existence proof, given in Appendix~\ref{appExtProof}, uses the notion of degree of divergence extended to $\mathcal{D}'(\Omega-O)$ as elaborated in the previous section. The rest of the lemma can be proved following the original result \cite{brunetti2000microlocal}.
\begin{lemma}
\label{ExistenceUniquenessExtensionLemma}
 Given $T_0 \in \mathcal{D}'(\Omega-O)$, with a finite degree of divergence, there exists an extension $T\in \mathcal{D}'(\Omega)$ such that $T|_{\Omega-O}=T_0$ and $\deg(T)=\deg({T_0})$. Moreover,
\begin{enumerate}[a.]
\item If $\deg{T_0} <0$ then the extension $T$ is unique.
\item  If $\deg{T_0} \geq 0$ then, given two extensions $T_1,T_2\in \mathcal{D}'(\Omega)$ which satisfy $T_1|_{\Omega-O}=T_2|_{\Omega-O}=T_0$ and $\deg(T_1)=\deg(T_2)=\deg{T_0}$, $T_1-T_2=\sum_{\alpha \in \mathbb{N}^n, |\alpha|\leq \deg(T_0)} T^{\alpha} \partial^\alpha \delta_O,$ where $T^{\alpha}\in \mathbb{R}.$
\end{enumerate}
\end{lemma}%
According to the lemma, a distribution $T\in \mathcal{D}'(\Omega)$ is uniquely characterized by $T|_{\Omega-O}$ if and only if $T$ has a negative degree of divergence. If $T$ has a non-negative degree of divergence, it can not be uniquely identified from $T|_{\Omega-O}$. The extent of non-uniqueness is as given in part \textit{b.} above.

\subsection{Singular fields}

We are interested in fields that are singular at a given point $O\in\Omega$. We say that a distribution $A\in \mathcal{A}(\Omega)\subset \mathcal{D}'(\Omega)$ if $\singsupp(A) \subset \{O\}$. The elements of $\mathcal{A}(\Omega)$ represent fields that are possibly singular at $O$ and smooth everywhere else. For $A\in \mathcal{A}(\Omega)$, $A|_{\Omega-O}$ can be represented in terms of a smooth function defined on $\Omega-O$, i.e., there exists a smooth map $a \in C^\infty (\Omega-O)$ such that $A|_{\Omega-O}(\psi)=\int_{\Omega-O} a \psi \dx$ for any $\psi \in \mathcal{D}(\Omega-O)$. 
We say that a distribution $E\in \mathcal{E}(\Omega)\subset \mathcal{D}'(\Omega)$ if $\supp(E) \subset \{O\}$.  For $E\in \mathcal{E}(\Omega)$ we have $E|_{\Omega-O}=0$, which implies that $\mathcal{E}(\Omega)$ is a subspace of $\mathcal{A}(\Omega)$ containing singular fields supported at $O$. 
For $A_1,A_2\in \mathcal{A}(\Omega)$, such that $A_1|_{\Omega-O}=A_2|_{\Omega-O}$, $(A_1-A_2)|_{\Omega-O}=0$ and hence $(A_1-A_2)\in\mathcal{E}(\Omega)$. Further, given two distributions $T\in \mathcal{D}'(\Omega)$ and and $E\in \mathcal{E}(\Omega)$, such that $T$ is smooth and $E\neq 0$, we have $(T+E)\in \mathcal{A}(\Omega)$ satisfying $\deg(T+E)\geq 0$.  According to the following lemma, any element of $\mathcal{E}(\Omega)$ can be written as a linear combination of $\delta_O$ and its derivatives. 
\begin{lemma}
\label{RepresentationELemma}
\textup{\cite[Theorem~3.2.1]{friedlander1998introduction}} For every $E\in \mathcal{E}(\Omega)$ we have the representation 
\begin{equation}
E=\sum_{\alpha \in \mathbb{N}^n, |\alpha|\leq \deg(E)} E^\alpha \partial^\alpha \delta_O, \label{representationE}
\end{equation}
with $E^\alpha \in \mathbb{R}$ given by $E^\alpha=E(w^\alpha)$, where $w^\alpha \in \mathcal{D}(\Omega)$ is such that, for any multi-index $\beta$, $\partial^\beta w^\alpha=(-1)^{|\alpha|}$ if $\alpha=\beta$ and $\partial^\beta w^\alpha=0$ if $\alpha \neq \beta$.
\end{lemma} 
An example of the compactly supported function $w^{\alpha}$, introduced in the above lemma, can be constructed by considering $\supp(w^{\alpha})\subset B_r$ with $w^{\alpha}(\boldsymbol{x})=(-1)^{|\alpha|}({x_1}^{\alpha_1} {x_2}^{\alpha_2}...{x_n}^{\alpha_n})/({\alpha_1 !\alpha_2 !...\alpha_n !})$ for all $\boldsymbol{x}\in B_{r/2}$.

For an open set $\Omega\subset \mathbb{R}^2$, given any $u\in \mathcal{D}'(\Omega)$, we have $\Curl (\nabla u)=0$. The converse, given in the following theorem, holds only when $\Omega$ is simply connected.
\begin{theorem}
\label{GeneralPoincareLemma}
\textup{\cite{mardare2008poincare}}
Let $\Omega \subset \mathbb{R}^2$ be simply connected, then for any $\boldsymbol{v} \in \mathcal{D}'(\Omega,\mathbb{R}^2)$ such that $\Curl \boldsymbol{v}=0$ there exists $u \in \mathcal{D}'(\Omega)$ such that $\nabla u =\boldsymbol{v}$.
\end{theorem}
Theorem~\ref{GeneralPoincareLemma} establishes the existence of a distribution $u$ such that $\nabla u=\boldsymbol{v}$ for any curl free vector valued distribution $\boldsymbol{v}$ on a simply connected open set. In the following lemma we establish more specific regularity results for $u$ when $\boldsymbol{v}$ is a singular field. For a connected open set $\Omega \subset \mathbb{R}^2$, any distribution $u\in \mathcal{D}'(\Omega)$ which satisfies $\nabla u =\boldsymbol{0}$ is a constant distribution, i.e., there exists a constant $c\in \mathbb{R}$ such that $u(\psi) = \int_\Omega c\psi \da$, where $\da$ denotes the area measure in $\mathbb{R}^2$, for any $\psi \in \mathcal{D}(\Omega)$. 
 \begin{lemma}
 \label{PLRegularity}
 Let $\Omega \subset \mathbb{R}^2$ be a simply connected open set. Then,
 \begin{enumerate}[a.]
 \item Given $\boldsymbol{A} \in \mathcal{A}(\Omega,\mathbb{R}^2)$, such that $\Curl \boldsymbol{A}=0$, there exists $u \in \mathcal{A}(\Omega)$ such that $\nabla u = \boldsymbol{A}$.
 \item Given $\boldsymbol{E} \in \mathcal{E}(\Omega,\mathbb{R}^2)$, such that $\Curl \boldsymbol{E}=0$,  there exists $u \in \mathcal{E}(\Omega)$ such that $\nabla u = \boldsymbol{E}$.
 \end{enumerate}

\begin{proof} \textit{a.} From Theorem~\ref{GeneralPoincareLemma} we have $u \in \mathcal{D}'(\Omega)$ such that $\nabla u = \boldsymbol{A}$. Thereafter, $\boldsymbol{A}|_{\Omega-O} = \left( \nabla u \right)|_{\Omega-O}=\nabla \left(u |_{\Omega-O} \right)$. The smoothness of $\boldsymbol{A}|_{\Omega-O}$ implies the smoothness of ${u}|_{\Omega-O}$, hence $u \in \mathcal{A}(\Omega)$.
\textit{b.} From Theorem~\ref{GeneralPoincareLemma} we have $u \in \mathcal{D}'(\Omega)$ such that $\nabla u = \boldsymbol{E}$. Then $\nabla \left(u |_{\Omega-O} \right) = \left( \nabla u \right)|_{\Omega-O}=0$.
Hence $u |_{\Omega-O}$ is equal to a constant distribution $c$ in $\Omega - O$. The distribution $u - c$ satisfies $\nabla (u - c)=\boldsymbol{E}$ and $\supp(u-c)=\{O\}$, which proves our assertion.
\end{proof}
\end{lemma}

The next lemma provides implications on the restricted smooth map, away from the point of singularity, given a vector valued distribution in $\mathcal{A}(\Omega,\mathbb{R}^2)$ which is curl or divergence free.
\begin{lemma}
\label{CurlDivRestrictionLemma}
Consider an open set $\Omega \subset \mathbb{R}^2$ and a distribution $\boldsymbol{A} \in \mathcal{A}(\Omega,\mathbb{R}^2)$. Then,
\begin{enumerate}[a.]
\item For $\Curl \boldsymbol{A}=0$ the restriction $\boldsymbol{A}|_{\Omega-O}$ satisfies
$\curl (\boldsymbol{A}|_{\Omega-O})=0$ and $\int_{\partial B_\epsilon} \left\langle\boldsymbol{A}|_{\Omega-O}, \boldsymbol{t}\right\rangle \dl=0$,
where $\boldsymbol{t}$ is the unit tangent to $\partial B_\epsilon$ and $\dl$ is the length measure in $\mathbb{R}^2$.
\item For $\Div \boldsymbol{A}=0$ the restriction $\boldsymbol{A}|_{\Omega-O}$ satisfies
$\div (\boldsymbol{A}|_{\Omega-O})=0$ and $\int_{\partial B_\epsilon} \left\langle\boldsymbol{A}|_{\Omega-O}, \boldsymbol{n}\right\rangle \dl=0$,
where $\boldsymbol{n}$ is the unit normal to $\partial B_\epsilon$.
\end{enumerate}
\begin{proof}
\textit{a.} Restricting $\Curl \boldsymbol{A}=0$ to $\Omega-O$ we obtain $\curl (\boldsymbol{A}|_{\Omega-O})=0$. Using Lemma~\ref{PLRegularity} and $\Curl \boldsymbol{A}|_{B_\epsilon}=0$ we have a distribution $u \in \mathcal{A}(B_\epsilon)$ such that $\nabla u = \boldsymbol{A}|_{B_\epsilon}$. Hence $\nabla \left(u |_{B_\epsilon-O} \right) = \boldsymbol{A}|_{B_\epsilon-O}$ in $B_\epsilon-O$, which immediately leads to the integral formula over $\partial B_\epsilon$.
\textit{b.} Noting that $\Div \boldsymbol{A}=0$ implies $\Curl (\boldsymbol{e}_3\times \boldsymbol{A})=0$, we can use the first part of the lemma to establish the required assertion. 
\end{proof}
\end{lemma}
In order to extend these results to the cases when $\Curl \boldsymbol{A}=E$ and $\Div \boldsymbol{A}=E$, where $E \in \mathcal{E}(\Omega)$, we first need the following lemma.
\begin{lemma}
\label{CurlDivDipoleSource}
Consider an open set $\Omega\subset \mathbb{R}^2$ and a distribution $E \in \mathcal{E}(\Omega)$ such that 
$ E^\alpha=0$, for $\alpha=(0,0)$, in the representation \eqref{representationE} for $E$ (with $n=2$). Then there exists $\boldsymbol{E}_1\in \mathcal{E}(\Omega,\mathbb{R}^2)$ and $\boldsymbol{E}_2\in \mathcal{E}(\Omega,\mathbb{R}^2)$ such that
$\Div \boldsymbol{E}_1=E$ and $\Curl \boldsymbol{E}_2=E$.
\begin{proof}
Given any multi index $\alpha=(\alpha_1,\alpha_2)$, such that $\alpha \neq (0,0)$, let $E_0= E^\alpha \partial^\alpha \delta_O.$ If $\alpha_1 \neq 0$ we take $\boldsymbol{E}_1= E^\alpha \partial^{\alpha'} \delta_O \boldsymbol{e}_1$ and $\boldsymbol{E}_2= E^\alpha \partial^{\alpha'} \delta_O \boldsymbol{e}_2$, where $\alpha'=(\alpha_1-1,\alpha_2)$. If $\alpha_1 = 0$ we take $\boldsymbol{E}_1= E^\alpha \partial^{\alpha'} \delta_O \boldsymbol{e}_2$ and $\boldsymbol{E}_2= -E^\alpha \partial^{\alpha'} \delta_O \boldsymbol{e}_1$, $\alpha'=(0,\alpha_2 -1)$. In either case, $\Div \boldsymbol{E}_1=E_0$ and $\Curl \boldsymbol{E}_2=E_0$.
\end{proof}
\end{lemma}
The next lemma obtains implications on the restricted smooth map, away from the point of singularity, given a vector valued distribution in $\mathcal{A}(\Omega,\mathbb{R}^2)$ whose curl or divergence belongs to $\mathcal{E}(\Omega)$.
\begin{lemma}
\label{CurlDivSingSupport}
Consider an open set $\Omega \subset \mathbb{R}^2$ and distributions $\boldsymbol{A} \in \mathcal{A}(\Omega,\mathbb{R}^2)$ and $E \in \mathcal{E}(\Omega)$, with the latter having a representation \eqref{representationE} with $n=2$. Then,
\begin{enumerate}[a.]
\item For $\Curl \boldsymbol{A}=E$ the restriction $\boldsymbol{A}|_{\Omega-O}$ satisfies
$\curl (\boldsymbol{A}|_{\Omega-O})=0$  and $\int_{\partial B_\epsilon} \left\langle\boldsymbol{A}|_{\Omega-O}, \boldsymbol{t}\right\rangle \dl=E^{(0,0)}$.
\item For $\Div \boldsymbol{A}=E$ the restriction $\boldsymbol{A}|_{\Omega-O}$ satisfies
$\div (\boldsymbol{A}|_{\Omega-O})=0$ and $\int_{\partial B_\epsilon} \left\langle\boldsymbol{A}|_{\Omega-O}, \boldsymbol{n}\right\rangle \dl=E^{(0,0)}$.
\end{enumerate}
\begin{proof}
\textit{a.} Using Lemma~\ref{CurlDivDipoleSource} we can construct $\boldsymbol{A}_1=\boldsymbol{E}+({E^{(0,0)}}/{2\pi r})\boldsymbol{e}_\theta$, where $\boldsymbol{E}\in \mathcal{E}(\Omega,\mathbb{R}^2)$, such that $\Curl \boldsymbol{A}_1=E.$ Hence $\Curl (\boldsymbol{A}-\boldsymbol{A}_1)=0$. An application of Lemma~\ref{CurlDivRestrictionLemma} establishes the required result.
\textit{b.} Using Lemma~\ref{CurlDivDipoleSource} we can construct $\boldsymbol{A}_1=\boldsymbol{E}+({E^{(0,0)}}/{2\pi r})\boldsymbol{e}_r$, where $\boldsymbol{E}\in \mathcal{E}(\Omega,\mathbb{R}^2)$, such that $\Div \boldsymbol{A}_1=E.$ Hence $\Div (\boldsymbol{A}-\boldsymbol{A}_1)=0$. An application of Lemma~\ref{CurlDivRestrictionLemma} establishes the required result.
\end{proof}
\end{lemma}
In the final lemma of this section, we establish conditions on a curl free restricted smooth map, away from the point of singularity, such that it yields a curl free extension over the whole domain.
\begin{lemma}
Given a smooth map $\boldsymbol{A}_0: \Omega-O \to \mathbb{R}^2$, with finite scaling degree, there exists an extension $\boldsymbol{A}\in \mathcal{A}(\Omega,\mathbb{R}^2)$ such that $\Curl \boldsymbol{A}={0}$ if and only if
$\curl \boldsymbol{A}_0={0}$ and 
$\int_{\partial B_\epsilon} \langle\boldsymbol{A}_0,\boldsymbol{t} \rangle \dl=0$.
\begin{proof}
The forward assertion has been established in Lemma~\ref{CurlDivRestrictionLemma}. On the other hand, for the given $\boldsymbol{A}_0$, there exists a distribution $\boldsymbol{A}_1 \in \mathcal{A}(\Omega,\mathbb{R}^2)$, such that $\boldsymbol{A}_1|_{\Omega-O}=\boldsymbol{A}_0$, satisfying $\curl \boldsymbol{A}_1|_{\Omega-O}={0}$ and 
$\int_{\partial B_\epsilon} \langle\boldsymbol{A}_1|_{\Omega-O},\boldsymbol{t} \rangle \dl=0$ (Lemma~\ref{ExistenceUniquenessExtensionLemma}). The former relation implies $\Curl \boldsymbol{A}_1 \in \mathcal{E}(\Omega)$. Then, as  a consequence of the latter combined with Lemma~\ref{CurlDivSingSupport}, $(\Curl \boldsymbol{A}_1)^{(0,0)} =0$ in the representation \eqref{representationE} for $\Curl \boldsymbol{A}_1$. Subsequently, in accordance with Lemma~\ref{CurlDivDipoleSource}, there exists $\boldsymbol{E} \in \mathcal{E}(\Omega,\mathbb{R}^2)$ such that $\Curl \boldsymbol{E}=\Curl \boldsymbol{A}_1$. The required distribution $\boldsymbol{A}\in \mathcal{A}(\Omega,\mathbb{R}^2)$, defined as $\boldsymbol{A}=\boldsymbol{A}_1 - \boldsymbol{E}$, satisfies $\boldsymbol{A}|_{\Omega-O}=\boldsymbol{A}_0$ and $\Curl \boldsymbol{A}=0$.
 \end{proof}
\end{lemma}

\subsection{An example of a singular field}

Consider an open set $\Omega \subset \mathbb{R}^2$. The functions $r$, $\cos \theta$, and $\sin \theta$ are smooth in $\Omega-O$ and locally integrable in $\Omega$ but not smooth in $\Omega$. Each of them can be used to define a distribution belonging to $\mathcal{A}(\Omega)$ with negative degree of divergence. In what follows, we construct a distribution $T\in \mathcal{A}(\Omega)$ such that $T|_{\Omega-O}$ is not locally integrable at $O$.
Let $G:\Omega-O \subset \mathbb{R}^2 \to \mathbb{R}$ be such that $G={g(\theta)}/{r^m}$, where $g$ is a bounded, smooth, and periodic non-trivial function of $\theta$. Since $g$ is periodic, with period $2\pi$, we can write $g(\theta)=\Sigma_{n=0}^{\infty} (c_n \cos(n\theta)+c'_n \sin(n\theta))$. Clearly, $G$ is not locally integrable at $O\in \Omega$ for $m\geq 2$. However, $G$ is smooth at every point in $\Omega-O$. Given $\phi \in \mathcal{D}(\Omega)$, consider the following functional associated with $G$:
\begin{equation}
T_G (\phi)=\lim_{\epsilon\to 0}\int_{\Omega-B_\epsilon} G \phi \da=\lim_{\epsilon\to 0}\int_{\Omega-B_\epsilon} \frac{g(\theta)}{r^m} \phi \da. 
\end{equation}
The functional $T_G$ will be well defined if the limit is well defined for arbitrary $\phi \in \mathcal{D} (\Omega)$. When $G$ is locally integrable at $O$, $T_G (\phi)=\int_{\Omega} G \phi \da$.
Since $\phi$ is smooth and compactly supported, for any natural number $k$ there exists a polynomial $P_k$ such that $|\phi(\boldsymbol{x})-P_k(\boldsymbol{x})| < K_k r^k$, where $K_k \in \mathbb{R}$ is finite, for all $\boldsymbol{x} \in B_r$ and some $r>0$.
We write
\begin{equation}
T_G (\phi)= \lim_{\epsilon\to 0} \left( \int_{\Omega-B_\epsilon} \frac{g(\theta)}{r^m} P_k \da + \int_{\Omega-B_\epsilon} \frac{g(\theta)}{r^m} (\phi-P_k) \da \right). 
\end{equation} 
For $k>m$ there always exists a polynomial $P_k$ such that the integrand in the second integral above is integrable and the limit exists. We are done if we establish the existence of the limit in the first term for arbitrary polynomials $P_k$. If the limit exists for arbitrary homogeneous polynomials it will exist for arbitrary polynomials. An arbitrary homogeneous polynomial of degree $n$ in $\mathbb{R}^2$ is of the form
$P_n(\boldsymbol{x})= r^n \langle \boldsymbol{a}_n, \boldsymbol{\mathit{E}}_n \rangle$,
where $\boldsymbol{\mathit{E}}_n=\boldsymbol{e}_r\otimes \cdots \text{(n times)} \cdots \otimes\boldsymbol{e}_r$ and $\boldsymbol{a}_n\in (\mathbb{R}^2)^n$. We evaluate
\begin{equation}
\lim_{\epsilon\to 0} \left( \int_{B_h-B_\epsilon} \frac{g(\theta)}{r^m} P_n \da \right)=\left(\lim_{\epsilon\to 0}\int_{\epsilon}^{h}{r^{(k-m+1)}} \text{d}r \right)  \left\langle \boldsymbol{a}_k,  \int_{0}^{2\pi} g(\theta) \boldsymbol{\mathit{E}}_k \text{d}\theta \right\rangle.
\end{equation}
On the right hand side of the expression, $\lim_{\epsilon\to 0}\int_{\epsilon}^{h}{r^{(k-m+1)}}\text{d}r$ exists only for $k>m-2$. For the limit to exist for arbitrary polynomials the necessary and sufficient condition is given by
$\int_{0}^{2\pi} g(\theta) \boldsymbol{\mathit{E}}_k \text{d}\theta = \boldsymbol{0}$ for all $k\leq m-2$.
Noting that $\int_{0}^{2\pi} \left(\cos n \theta \right) \boldsymbol{\mathit{E}}_k \text{d}\theta =\boldsymbol{0}$
and 
$\int_{0}^{2\pi} \left( \sin n \theta \right)  \boldsymbol{\mathit{E}}_k \text{d}\theta =\boldsymbol{0}$
for all $n>k$, we obtain the conditions for the limit to exist as $c_n=c'_n=0$ for all $n\leq m-2$.

\section{Stress equilibrium} \label{stressequibsec}

\subsection{A generalized equilibrium condition}
In order to allow both the stress field $\boldsymbol{\sigma}$ and the body force field $\boldsymbol{B}$ to develop singularities at isolated points in a domain $\Omega \subset \mathbb{R}^2$ we consider them as distributions $\boldsymbol{\sigma}\in \mathcal{D}'(\Omega,\Sym)$ and $\boldsymbol{B}\in \mathcal{D}'(\Omega,\mathbb{R}^2)$. Moreover, in the absence of inertial forces, we postulate the fields to satisfy an equilibrium condition given in a distributional form as
\begin{equation}
\label{DistributionalBalanceLaw}
\Div \boldsymbol{\sigma} + \boldsymbol{B}=\boldsymbol{0}.
\end{equation}
If the fields are smooth over $\Omega$ then this reduces to the classical form of a pointwise balance law. Before we move on towards studying implications of \eqref{DistributionalBalanceLaw} for singular fields, we would like to emphasize the generality afforded by the given form of the equilibrium condition. Conventionally, for weakly regular stress fields (such as those not necessarily continuous over $\Omega$), the equilibrium is postulated using a Cauchy flux map $\mathcal{F}$ such that 
\begin{equation}
\label{CauchyFluxasLimit}
\mathcal{F}(\mathcal{S})=\lim_{\rho \to 0} \int_{\mathcal{S}} \boldsymbol{\sigma}_\rho \boldsymbol{n} \da
\end{equation}
for any smooth oriented surface $\mathcal{S}$ in $\Omega$ with unit normal $\boldsymbol{n}$, where $\boldsymbol{\sigma}_\rho$ is a sequence of smooth fields such that $\boldsymbol{\sigma}_\rho \to \boldsymbol{\sigma}$ in the sense of distributions \cite{vsilhavy1987existence, vsilhavy2008cauchy}. For a smooth stress field, $\mathcal{F}(\mathcal{S})=\int_{\mathcal{S}} \boldsymbol{\sigma}\boldsymbol{n} \da$.
The map $\mathcal{F}(\mathcal{S})$ represents the contact force transmitted across the surface $\mathcal{S}$. Given a map $\mathcal{F}$ and a smooth body force field $\boldsymbol{b}:\Omega\to \mathbb{R}^2$, the equilibrium condition is given by
\begin{equation}
\label{BalanceLaw}
\mathcal{F}(\partial \mathcal{P})+\int_\mathcal{P} \boldsymbol{b} \da=\boldsymbol{0},
\end{equation}
where $\mathcal{P} \subset \Omega$ is an arbitrary open subset of $\Omega$ with smooth boundary $\partial \mathcal{P}$. For a smooth stress field this is equivalent to the point wise condition $\div \boldsymbol{\sigma} + \boldsymbol{b}=\boldsymbol{0}$ in $\Omega$. The symmetric nature of the smooth stress field is a consequence of the angular momentum balance.
For a general stress field $\boldsymbol{\sigma}\in \mathcal{D}'(\Omega,\Sym)$, the limit in \eqref{CauchyFluxasLimit} does not exist for all smooth surfaces in $\Omega$ \cite{podio2006concentrated}. For instance, if $\boldsymbol{\sigma}=\delta_O \boldsymbol{I}$ then $\mathcal{F}$ is well defined only for surfaces $\mathcal{S}$ such that $O \notin \mathcal{S}$. In such a situation, the stress field $\boldsymbol{\sigma}\in \mathcal{D}'(\Omega,\Sym)$ can not be interpreted in terms of the Cauchy flux map and \eqref{BalanceLaw} can no longer be used as the general equilibrium condition since $\mathcal{F}(\partial \mathcal{P})$ is not defined for arbitrary $\mathcal{P}$.
We can think of the generalization of the equilibrium condition, as given in \eqref{DistributionalBalanceLaw}, in the following limiting sense. For any field $\boldsymbol{\sigma}\in \mathcal{D}'(\Omega,\Sym)$ there exists a sequence of smooth maps $\boldsymbol{\sigma}_\rho$ such that $\boldsymbol{\sigma}_\rho \to \boldsymbol{\sigma}$ as $\rho \to 0$ \cite[Section~5.2]{friedlander1998introduction}. Any distributional stress field can therefore be interpreted as the limit of a sequence of smooth stress fields. The distributional stress field is said to be in equilibrium if it is the limit of a sequence of smooth equilibrated stress fields. The body force field $\boldsymbol{B}\in \mathcal{D}'(\Omega,\mathbb{R}^2)$ is then the limiting value of the corresponding sequence of smooth body force fields $\boldsymbol{B}_\rho$. The equilibrium condition \eqref{DistributionalBalanceLaw} follows immediately as the limit of the conditions $\div \boldsymbol{\sigma}_\rho + \boldsymbol{B}_\rho=0$ in $\Omega$ as $\rho \to 0$. 

\subsection{Stress fields with point singularity}
The restriction $\boldsymbol{\sigma}|_{\Omega-O}$ of a singular stress field $\boldsymbol{\sigma}\in \mathcal{A}(\Omega,\Sym)$ is smooth. The Cauchy flux is then well defined for any surface contained in $\Omega-O$ but not necessarily so for surfaces $\mathcal{S}$ such that $O \in \mathcal{S}.$ If $\boldsymbol{\sigma}$ is in equilibrium with a general body force field $\boldsymbol{B}\in \mathcal{D}'(\Omega,\mathbb{R}^2)$ then $\boldsymbol{B}\in \mathcal{A}(\Omega,\mathbb{R}^2)$.  In this section, we obtain both the implications and the local equivalent relations of the equilibrium condition \eqref{DistributionalBalanceLaw}. We emphasize whenever the local relations can be written completely in terms of the smooth restriction $\boldsymbol{\sigma}|_{\Omega-O}$. We begin with
\begin{lemma}
\label{stressequiv1}
Let $\Omega \subset \mathbb{R}^2$. Consider a singular stress field $\boldsymbol{\sigma}\in \mathcal{A}(\Omega,\Sym)$ and a singular body force field $\boldsymbol{B}\in \mathcal{A}(\Omega,\mathbb{R}^2)$. Then,
\begin{enumerate}[a.]
\item The equilibrium condition \eqref{DistributionalBalanceLaw}
 implies the pointwise condition $\div(\boldsymbol{\sigma}|_{\Omega-O})+\boldsymbol{B}|_{\Omega-O}=\boldsymbol{0}$.
\item If $\deg(\boldsymbol{\sigma}) < -1$ and $\deg (\boldsymbol{B}) < 0$ then the equilibrium condition \eqref{DistributionalBalanceLaw} is equivalent to $\div(\boldsymbol{\sigma}|_{\Omega-O})+\boldsymbol{B}|_{\Omega-O}=\boldsymbol{0}$.
\end{enumerate}
\begin{proof}
\textit{a.} The result is obtained by restricting the general equilibrium condition \eqref{DistributionalBalanceLaw} to $\Omega-O$. 
\textit{b.} With the given degrees of divergence for $\boldsymbol{\sigma}$ and $\boldsymbol{B}$, $\deg(\Div \boldsymbol{\sigma} + \boldsymbol{B})<0$. Hence, by Lemma~\ref{ExistenceUniquenessExtensionLemma}, the unique extension of $\div(\boldsymbol{\sigma}|_{\Omega-O})+\boldsymbol{B}|_{\Omega-O}$ is $\Div \boldsymbol{\sigma} + \boldsymbol{B}$.
\end{proof}
\end{lemma}

We give two examples, one where Lemma~\ref{stressequiv1}\textit{b.} holds and one where it does not. Consider an integrable stress field $\boldsymbol{\sigma} \in \mathcal{D}'(\Omega,\Sym)$ such that $\Div \boldsymbol{\sigma}$ is also integrable. Therefore, there exist integrable functions $\boldsymbol{f}:\Omega \to \Sym$ and $\boldsymbol{g}:\Omega\to \mathbb{R}^2$ such that $\boldsymbol{\sigma}(\boldsymbol{\phi})=\int_\Omega \langle \boldsymbol{f},\boldsymbol{\phi} \rangle \da$, for all $\boldsymbol{\phi}\in \mathcal{D}(\Omega,\Sym)$, and $\Div \boldsymbol{\sigma}( \boldsymbol{\psi})=\int_\Omega \langle \boldsymbol{g},\boldsymbol{\psi} \rangle \da$, for all $\boldsymbol{\psi}\in \mathcal{D}(\Omega,\mathbb{R}^2)$. If the body force field $\boldsymbol{B}\in \mathcal{D}'(\Omega,\mathbb{R}^2)$ is integrable then $(\Div \boldsymbol{\sigma}+\boldsymbol{B})$ is also integrable. Assuming $\div(\boldsymbol{\sigma}|_{\Omega-O})+\boldsymbol{B}|_{\Omega-O}=\boldsymbol{0}$ we can subsequently conclude that  $(\Div \boldsymbol{\sigma}+\boldsymbol{B}) =\boldsymbol{0}$. Furthermore, if the stress field is singular, i.e.,  $\boldsymbol{\sigma} \in \mathcal{A}(\Omega,\Sym)$, then there exists a smooth map $\boldsymbol{\sigma}_0:\Omega-O\to \Sym$, such that  $\boldsymbol{\sigma}(\boldsymbol{\phi})=\int_\Omega \langle \boldsymbol{\sigma}_0,\boldsymbol{\phi} \rangle \da$, for all $\boldsymbol{\phi}\in \mathcal{D}(\Omega,\Sym)$, and $\Div\boldsymbol{\sigma}( \boldsymbol{\psi})=\int_\Omega \langle \div \boldsymbol{\sigma}_0,\boldsymbol{\psi} \rangle \da$, for all $\boldsymbol{\psi}\in \mathcal{D}(\Omega,\mathbb{R}^2)$. The stress equilibrium condition can therefore be completely given in terms of $\boldsymbol{\sigma}_0$.
On the other hand, consider a stress field of the form $\boldsymbol{\sigma}(\boldsymbol{\phi})=\int_{\Omega}({\cos \theta}/{r})\left\langle  \boldsymbol{\phi} ,\boldsymbol{e}_r \otimes \boldsymbol{e}_r\right\rangle \da$, for an arbitrary $\boldsymbol{\phi}\in \mathcal{D}(\Omega,\Sym)$, with $\boldsymbol{B}=\boldsymbol{0}$. Such a $\boldsymbol{\sigma}$ is integrable but $\Div \boldsymbol{\sigma}$ is not integrable. In this case, $(\Div \boldsymbol{\sigma}+\boldsymbol{B})|_{\Omega-O}=\boldsymbol{0}$ but $\Div \boldsymbol{\sigma}+\boldsymbol{B} \neq \boldsymbol{0}$. For the considered $\boldsymbol{\sigma}$, $\deg(\boldsymbol{\sigma})= -1$ which violates the assumption required for Lemma~\ref{stressequiv1}\textit{b.} 

In the next lemma, we discuss the equilibrium of a singular stress field with a point supported body force field.
\begin{lemma}
\label{stressequiv2}
Let $\Omega \subset \mathbb{R}^2$. Consider a singular stress field $\boldsymbol{\sigma}\in \mathcal{A}(\Omega,\Sym)$ and body force field $\boldsymbol{B}\in \mathcal{E}(\Omega,\mathbb{R}^2)$ with the representation $\boldsymbol{B}=\sum_{\alpha \in \mathbb{N}^2, |\alpha|\leq \deg(\boldsymbol{B})} \boldsymbol{b}^\alpha \partial^\alpha \delta_O$, where $\boldsymbol{b}^\alpha \in \mathbb{R}^2$. Then,
\begin{enumerate}[a.]
\item 
The equilibrium condition \eqref{DistributionalBalanceLaw}
 implies
\begin{subequations}\label{restrictedBalance11}%
\begin{align}
&\div(\boldsymbol{\sigma}|_{\Omega-O})=\boldsymbol{0}, \label{restrictedBalance1}\\
&\int_{\partial B_\epsilon}(\boldsymbol{\sigma}|_{\Omega-O} \boldsymbol{n}) \dl =-\boldsymbol{b}^{(0,0)},~\text{and} \label{restrictedBalance2} \\
& \deg(\boldsymbol{B})\leq \deg(\boldsymbol{\sigma})+1. \label{restrictedBalance3}
\end{align}
\end{subequations}
\item For $\deg(\boldsymbol{\sigma}) < 0$ the equilibrium condition \eqref{DistributionalBalanceLaw} is equivalent to Equations~\eqref{restrictedBalance11}.
\end{enumerate}
\begin{proof}
\textit{a.}  Equations~\eqref{restrictedBalance1} and \eqref{restrictedBalance2} follow from Lemma~\ref{CurlDivSingSupport} whereas Equation~\eqref{restrictedBalance3} follows from point~v. given at the end of Section~\ref{sddeg}.
\textit{b.}  Equation~\eqref{restrictedBalance3} with $\deg(\boldsymbol{\sigma}) < 0$ requires that $\deg(\boldsymbol{B})<1$. Then, necessarily $\boldsymbol{B}=\boldsymbol{b}^{(0,0)}\delta_O$. On the other hand, \eqref{restrictedBalance1} and $\deg(\Div \boldsymbol{\sigma})<1$, in conjunction with Lemma~\ref{RepresentationELemma}, imply that $\Div \boldsymbol{\sigma}=\boldsymbol{a} \delta_O$, where $\boldsymbol{a}\in \mathbb{R}^2$.  That $\boldsymbol{a}=-\boldsymbol{b}^{(0,0)}$, follows from \eqref{restrictedBalance2} and Lemma~\ref{CurlDivSingSupport}.
\end{proof}
\end{lemma}

Consequent to the preceding lemma, we provide several examples which illustrate sufficiency (or insufficiency) of Equations~\eqref{restrictedBalance11} in enforcing the stress equilibrium in $\Omega$. Consider a stress field $\boldsymbol{\sigma}\in \mathcal{A}(\Omega,\Sym)$ such that $\deg(\boldsymbol{\sigma})=\deg(\boldsymbol{\sigma}|_{\Omega-O})$ and $\boldsymbol{\sigma}|_{\Omega-O}$  behaves like $r^k$. Then $\deg(\boldsymbol{\sigma})=\deg(\boldsymbol{\sigma}|_{\Omega-O})=-k-2$. According to Lemma~\ref{stressequiv2}\textit{b.},   Equations~\eqref{restrictedBalance1} and \eqref{restrictedBalance2} are sufficient to enforce the equilibrium as long as $k > -2$. For instance, let $\boldsymbol{\sigma}$ be such that $\boldsymbol{\sigma}(\boldsymbol{\phi})=\int_{\Omega}({\cos \theta}/{r})\left\langle  \boldsymbol{\phi} ,\boldsymbol{e}_r\otimes \boldsymbol{e}_r \right\rangle \da$, for any $\boldsymbol{\phi}\in \mathcal{D}(\Omega,\Lin)$, and let $\boldsymbol{B}=-\pi\boldsymbol{e}_1\delta_O$. Therefore, $k=-1$, $\deg(\boldsymbol{\sigma})=-1$, and $\deg(\boldsymbol{B})=0$. Further, we calculate $\div (\boldsymbol{\sigma}|_{\Omega-O})=\boldsymbol{0}$ and $\int_{\partial B_\epsilon}(\boldsymbol{\sigma}|_{\Omega-O} \boldsymbol{n}) \dl=\pi \boldsymbol{e}_1$. These two relations are sufficient to ensure that $\Div \boldsymbol{\sigma} -\pi\boldsymbol{e}_1\delta_O = \boldsymbol{0}$ holds true.  On the other hand, if we consider a stress field such that $\boldsymbol{\sigma}(\boldsymbol{\phi})=\lim_{\epsilon \to 0} \int_{\Omega-B_\epsilon} ({1}/{\pi r^2}) \left\langle  \boldsymbol{\phi} ,-\boldsymbol{e}_r\otimes \boldsymbol{e}_r+\boldsymbol{e}_\theta\otimes \boldsymbol{e}_\theta \right\rangle \da$, for any $\boldsymbol{\phi}\in \mathcal{D}(\Omega,\Lin)$, and a body force field $\boldsymbol{B}=\boldsymbol{0}$, then
 $k=-2$ and $\deg(\boldsymbol{\sigma})=0$. Both \eqref{restrictedBalance1} and \eqref{restrictedBalance2} are trivially satisfied but $\Div \boldsymbol{\sigma} = {\nabla \delta_O}$, which is inconsistent with the equilibrium condition since $\boldsymbol{B}=\boldsymbol{0}$. Therefore, although Equations~\eqref{restrictedBalance1} and \eqref{restrictedBalance2} are satisfied, the given stress field is not equilibrated. These equations are therefore insufficient in determining whether the given stress field is in equilibrium or not. A similar conclusion is realized when we consider point supported singular stress fields $\boldsymbol{\sigma}\in \mathcal{E}(\Omega,\Sym)$. For such fields $\deg(\boldsymbol{\sigma})\geq 0$ and Equations~\eqref{restrictedBalance1} and \eqref{restrictedBalance2} are satisfied trivially as long as we consider body force fields with $\boldsymbol{b}^{(0,0)} =\boldsymbol{0}$. This however has no bearing on whether the condition $\Div \boldsymbol{\sigma}+\boldsymbol{B}=\boldsymbol{0}$ holds true or not. Finally, consider $\boldsymbol{B}\in \mathcal{E}(\Omega,\mathbb{R}^2)$ such that $\deg (\boldsymbol{B})\geq1$ without any restriction on $\boldsymbol{\sigma}$. An example of such a body force field is the concentrated force dipole, where $\boldsymbol{B}=(\boldsymbol{e}_1 \otimes \boldsymbol{e}_1) {\nabla \delta_O}$ and $\deg (\boldsymbol{B})=1$. In such a case \eqref{restrictedBalance3} requires $\deg(\boldsymbol{\sigma})\geq 0$. As already shown through a counter example above, Equations~\eqref{restrictedBalance1} and \eqref{restrictedBalance2} are again not sufficient for determining the veracity of stress equilibrium.

It may seem from the above discussion that there could be additional conditions on $\boldsymbol{\sigma}|_{\Omega-O}$, besides \eqref{restrictedBalance1} and \eqref{restrictedBalance2}, which need to be enforced in order to establish the equilibrium of stress. In fact, in general, this is not so and one would need the complete stress field $\boldsymbol{\sigma}\in \mathcal{A}(\Omega,\Sym)$ for the consideration of equilibrium. This claim is elaborated in the following remark.
\begin{remark} \label{stressoutinsuf}
 Given a body force field $\boldsymbol{B}\in \mathcal{D}'(\Omega,\mathbb{R}^2)$, let $\boldsymbol{\sigma}_0:\Omega-O \to \Sym$ be a smooth map of finite scaling degree such that $\div\boldsymbol{\sigma}_0+\boldsymbol{B}|_{\Omega-O}=\boldsymbol{0}.$ For any extension $\boldsymbol{\sigma}\in \mathcal{A}(\Omega,\Sym)$ of $\boldsymbol{\sigma}_0$ we have $(\Div\boldsymbol{\sigma}+\boldsymbol{B}) \in \mathcal{E}(\Omega,\mathbb{R}^3)$. We assume that there exists an extension $\boldsymbol{\sigma}_1 \in \mathcal{A}(\Omega,\Sym)$ which is in equilibrium with $\boldsymbol{B}$, i.e., $\Div \boldsymbol{\sigma}_1+\boldsymbol{B}=\boldsymbol{0}$. For any $\boldsymbol{\sigma}_2\in \mathcal{E}(\Omega,\Sym)$, such that $\Div \boldsymbol{\sigma}_2 \neq \boldsymbol{0}$, the field $\boldsymbol{\sigma}_3=\boldsymbol{\sigma}_1 +\boldsymbol{\sigma}_2$ is an extension of $\boldsymbol{\sigma}_0$ which is not in equilibrium with $\boldsymbol{B}$, i.e.,
 $\Div \boldsymbol{\sigma}_3+\boldsymbol{B}\neq \boldsymbol{0}$.
Hence no equilibrium condition on $\boldsymbol{\sigma}_0$ would, in general, guarantee equilibrium for an arbitrary extension of $\boldsymbol{\sigma}_0$. 
\end{remark}

In the following lemma we obtain the complete set of local conditions on a singular stress field which are equivalent to the equilibrium condition \eqref{DistributionalBalanceLaw}.
 
 \begin{lemma} \label{fullocalstressequib}
Let $\Omega \subset \mathbb{R}^2$. For a singular stress field $\boldsymbol{\sigma}\in \mathcal{A}(\Omega,\Sym)$ and a singular body force field $\boldsymbol{B}\in \mathcal{A}(\Omega,\Sym)$, the equilibrium condition \eqref{DistributionalBalanceLaw} is equivalent to
\begin{subequations}
\begin{align}
&\div(\boldsymbol{\sigma}|_{\Omega-O})+\boldsymbol{B}|_{\Omega-O}=\boldsymbol{0}, \label{restrictedBalance_1} \\
&- \boldsymbol{\sigma}(\nabla (w^\alpha \boldsymbol{e}_1))+ \boldsymbol{B}( w^\alpha \boldsymbol{e}_1)=0,  \label{restrictedBalance_2}\\
&- \boldsymbol{\sigma}(\nabla (w^\alpha \boldsymbol{e}_2))+ \boldsymbol{B}( w^\alpha \boldsymbol{e}_2)=0,~\text{and}  \label{restrictedBalance_3}\\
 &\deg(\boldsymbol{B})\leq \deg(\boldsymbol{\sigma})+1  \label{restrictedBalance_4}
 \end{align}
 \label{restrictedBalance22}
 \end{subequations}
 for all $|\alpha|\leq\deg(\boldsymbol{\sigma})+1$, where $\alpha$ is a multi-index and $w^\alpha \in \mathcal{D}(\Omega)$ is as introduced in Lemma~\ref{RepresentationELemma}.
 
\begin{proof}
Given \eqref{DistributionalBalanceLaw}, Equation~\eqref{restrictedBalance_1} follows by restricting it to $\Omega-O$, \eqref{restrictedBalance_4} follows from point~v. given at the end of Section~\ref{sddeg}, and Equations~\eqref{restrictedBalance_2} and \eqref{restrictedBalance_3} follow after using $w^\alpha \boldsymbol{e}_1$ and $w^\alpha \boldsymbol{e}_2$ as test functions, respectively, in the distributional form of \eqref{DistributionalBalanceLaw}. On the other hand, Equation~\eqref{restrictedBalance_1} implies that $(\Div \boldsymbol{\sigma}+\boldsymbol{B})\in \mathcal{E}(\Omega,\mathbb{R}^2)$. Due to \eqref{restrictedBalance_4}, and the fact that $\deg(\Div \boldsymbol{\sigma}) \leq (\deg(\boldsymbol{\sigma})+1)$, we have $\deg(\Div \boldsymbol{\sigma}+\boldsymbol{B})\leq (\deg(\boldsymbol{\sigma})+1)$.  Equations~\eqref{restrictedBalance_2} and \eqref{restrictedBalance_3}, in conjunction with  Lemma~\ref{RepresentationELemma}, subsequently yield \eqref{DistributionalBalanceLaw}.
\end{proof} 
 \end{lemma}
 
Equations~\eqref{restrictedBalance_2} and \eqref{restrictedBalance_3} are local at $O$ in the following sense. For any two singular fields $\boldsymbol{\sigma}_1,\boldsymbol{\sigma}_2\in \mathcal{A}(\Omega,\Sym)$ such that $\boldsymbol{\sigma}_1|_\omega=\boldsymbol{\sigma}_2|_\omega$, where $\omega \subset \Omega$ is an arbitrary open set with $O\in \omega$, \eqref{restrictedBalance_2} and \eqref{restrictedBalance_3} are equivalently written in terms of either $\boldsymbol{\sigma}_1$ or $\boldsymbol{\sigma}_2$. If $k$ is the greatest integer smaller (or equal) than $\deg(\boldsymbol{\sigma})$, with $\deg(\boldsymbol{\sigma})\geq -1$, then the total number of conditions represented by \eqref{restrictedBalance_2} and \eqref{restrictedBalance_3}  is given as $(2^{k+3}-2)$. 
For instance, when $-1\leq\deg(\boldsymbol{\sigma})<0,$ these yield two scalar conditions. The two conditions are equivalent to the loop integral condition \eqref{restrictedBalance2} when $\boldsymbol{B}|_{\Omega-O}=\boldsymbol{0}$. As an example to illustrate the lemma, consider  $\boldsymbol{\sigma}(\boldsymbol{\phi})=\lim_{\epsilon \to 0} \int_{\Omega-B_\epsilon} ({1}/{\pi r^2}) \left\langle  \boldsymbol{\phi} ,-\boldsymbol{e}_r\otimes \boldsymbol{e}_r+\boldsymbol{e}_\theta\otimes \boldsymbol{e}_\theta \right\rangle \da$, for any $\boldsymbol{\phi}\in \mathcal{D}(\Omega,\Lin)$, and $\boldsymbol{B}=\boldsymbol{0} $. Clearly, Equation~\eqref{restrictedBalance_1} is satisfied but \eqref{restrictedBalance_2} and \eqref{restrictedBalance_3} hold true only for $|\alpha|=0$ but not for $|\alpha|=1$. The given fields therefore are not in equilibrium. However, the same stress field with body force $\boldsymbol{B}=-\nabla \delta_O$ satisfies all the conditions in \eqref{restrictedBalance22} thereby validating the stress equilibrium condition.

\section{Strain compatibility and incompatibility} \label{straincompincomp}

\subsection{A general strain compatibility condition}
In order to allow the strain field $\boldsymbol{E}$ to have singularities at isolated points in $\Omega \subset \mathbb{R}^2$ we consider it as a distribution in $\mathcal{D}'(\Omega,\Sym)$. A strain field $\boldsymbol{E}\in \mathcal{D}'(\Omega,\Sym)$ is said to be compatible if and only if there exists $\boldsymbol{u}\in \mathcal{D}'(\Omega,\mathbb{R}^2)$ such that $\boldsymbol{E}=(1/2)\left({\nabla \boldsymbol{u} + ( \nabla \boldsymbol{u} )^T}\right)$. The field $\boldsymbol{u}$ is the displacement associated with the strain $\boldsymbol{E}$. The necessary and sufficient conditions for a strain to be compatible over a simple connected domain are given in the lemma below.
\begin{lemma} \textup{\cite{pandey2020topological}} \label{strcompl1}
For a simply connected region $\Omega\subset \mathbb{R}^2$ and a strain field $\boldsymbol{E} \in \mathcal{D}'(\Omega,\Sym)$, there exists $\boldsymbol{u} \in \mathcal{D}'(\Omega,\mathbb{R}^2)$ such that $\boldsymbol{E}=(1/2)\left({\nabla \boldsymbol{u} + ( \nabla \boldsymbol{u} )^T}\right)$ if and only if 
\begin{equation}
\Curl\Curl \boldsymbol{E}={0}. \label{compatcond}
\end{equation} 
\end{lemma}
Moreover, given two displacement fields $\boldsymbol{u}_1,\boldsymbol{u}_2 \in \mathcal{D}'(\Omega,\mathbb{R}^2)$, such that ${\nabla \boldsymbol{u}_1 + \left( \nabla \boldsymbol{u}_1 \right)^T}= {\nabla \boldsymbol{u}_2 + \left( \nabla \boldsymbol{u}_2 \right)^T}$, there exist constants $\boldsymbol{c}_0\in \mathbb{R}^2$ and $c_1 \in \mathbb{R}$ such that $\boldsymbol{u}_2=\boldsymbol{u}_1+\boldsymbol{c}_0+c_1 (\boldsymbol{e}_3\times \boldsymbol{x}).$

\subsection{Compatibility of singular strain fields}
\label{compsingstr}

We begin by establishing stronger versions of Lemma~\ref{strcompl1} when $\boldsymbol{E} \in \mathcal{A}(\Omega,\Sym)$ or $\boldsymbol{E} \in \mathcal{E}(\Omega,\Sym)$.
\begin{lemma} \label{strcompl2} For a simply connected region $\Omega\subset \mathbb{R}^2$ and a strain field $\boldsymbol{E}\in \mathcal{D}'(\Omega,\Sym)$ which satisfies $\Curl \Curl \boldsymbol{E}=0$,
\begin{enumerate}[a.]
\item  If $\boldsymbol{E} \in \mathcal{A}(\Omega,\Sym)$ then there exists $\boldsymbol{u} \in \mathcal{A}(\Omega,\mathbb{R}^2)$ such that $\boldsymbol{E}=(1/2)\left({\nabla \boldsymbol{u} + ( \nabla \boldsymbol{u} )^T}\right)$.
\item If $\boldsymbol{E} \in \mathcal{E}(\Omega,\Sym)$ then there exists $\boldsymbol{u} \in \mathcal{E}(\Omega,\mathbb{R}^2)$ such that $\boldsymbol{E}=(1/2)\left({\nabla \boldsymbol{u} + ( \nabla \boldsymbol{u} )^T}\right)$.
\end{enumerate}
\begin{proof}
\textit{a.} We have $\boldsymbol{u} \in \mathcal{D}'(\Omega,\mathbb{R}^2)$ such that $\boldsymbol{E}=(1/2)\left({\nabla \boldsymbol{u} + ( \nabla \boldsymbol{u} )^T}\right)$. Restricting to $\Omega -O$ we write $\boldsymbol{E}|_{\Omega-O} = (1/2)\left({\nabla (\boldsymbol{u} |_{\Omega-O})+(\nabla (\boldsymbol{u} |_{\Omega-O}))^T } \right)$.
The smoothness of $\boldsymbol{E}|_{\Omega-O}$ implies the smoothness of ${\boldsymbol{u}}\arrowvert _{\Omega-O}$, hence $\boldsymbol{u} \in \mathcal{A}(\Omega, \mathbb{R}^2)$.
\textit{b.} We have $\boldsymbol{u}_1 \in \mathcal{D}'(\Omega,\mathbb{R}^2)$ such that $\boldsymbol{E}=(1/2)\left({\nabla \boldsymbol{u}_1 + ( \nabla \boldsymbol{u}_1 )^T}\right)$. Restricting to $\Omega -O$ we obtain
$(1/2) \left({\nabla \boldsymbol{u}_1|_{\Omega-O}+( \nabla \boldsymbol{u}_1|_{\Omega-O})^T}\right)  = 0$.
Hence $\boldsymbol{u}_1 |_{\Omega-O}=\boldsymbol{c}_0 +c_1 \boldsymbol{e}_3\times \boldsymbol{x}$, where $\boldsymbol{c}_0\in \mathbb{R}^2$ and $c_1 \in \mathbb{R}$ are constants. The field $\boldsymbol{u}=\boldsymbol{u}_1 - \boldsymbol{c}_0-c_1 \boldsymbol{e}_3\times \boldsymbol{x}$ satisfies $\boldsymbol{E}=(1/2)\left({\nabla \boldsymbol{u} + ( \nabla \boldsymbol{u} )^T}\right)$ with $\supp(\boldsymbol{u})=\{O\}$.
\end{proof}
\end{lemma}
In the next lemma we obtain necessary and sufficient conditions on a smooth map $\boldsymbol{E}_0: \Omega-O \to \Sym$ such that it admits a compatible extension in $\Omega$.
\begin{lemma}
\label{strcompl3} 
Given a simply connected region $\Omega\subset \mathbb{R}^2$ and a smooth map $\boldsymbol{E}_0: \Omega-O \to \Sym$, with finite scaling degree, there exists an extension $\boldsymbol{E}\in \mathcal{A}(\Omega,\Sym)$ such that $\Curl \Curl \boldsymbol{E}={0}$ if and only if
\begin{subequations}
\begin{align}
&\curl \curl \boldsymbol{E}_0={0} ~\text{and}\label{RestrictionCurlCurl1} \\
& \int_{\partial B_\epsilon} \{ \boldsymbol{E}_0(\boldsymbol{y})  + (\boldsymbol{y}-\boldsymbol{x})\times \curl (\boldsymbol{E}_0(\boldsymbol{y})) \} \dy =\boldsymbol{0} \label{RestrictionCurlCurl2}
\end{align}
\end{subequations}
for all $\boldsymbol{x} \in \Omega-O$.
\begin{proof}
For $\boldsymbol{E}\in \mathcal{A}(\Omega,\Sym)$ satisfying $\Curl \Curl \boldsymbol{E}={0}$, Lemma~\ref{strcompl2} implies that there exists $\boldsymbol{u} \in \mathcal{A}(\Omega,\mathbb{R}^2)$ such that $\boldsymbol{E}=(1/2)\left({\nabla \boldsymbol{u} + ( \nabla \boldsymbol{u} )^T}\right)$.  Accordingly, $\boldsymbol{E}_0=(1/2)\left({\nabla (\boldsymbol{u}|_{\Omega-O})+\nabla (\boldsymbol{u}|_{\Omega-O})^T} \right)$ which immediately leads to \eqref{RestrictionCurlCurl1} and \eqref{RestrictionCurlCurl2}. On the other hand,
given $\boldsymbol{E}_0$, which satisfies conditions from the lemma, there exists a smooth map $\boldsymbol{u}_0: \Omega-O \to \mathbb{R}^2$, of finite scaling degree, such that $\boldsymbol{E}_0=(1/2)\left({\nabla \boldsymbol{u}_0 + ( \nabla \boldsymbol{u}_0 )^T}\right)$. Following Lemma~\ref{ExistenceUniquenessExtensionLemma}, we have $\boldsymbol{u} \in \mathcal{A}(\Omega, \mathbb{R}^2)$ such that $\boldsymbol{u}|_{\Omega-O}=\boldsymbol{u}_0$. A distribution  $\boldsymbol{E}\in \mathcal{A}(\Omega,\Sym)$, defined as $\boldsymbol{E}=(1/2)\left({\nabla \boldsymbol{u} + ( \nabla \boldsymbol{u} )^T}\right)$ satisfies both $\boldsymbol{E}|_{\Omega-O}=\boldsymbol{E}_0$ and $\Curl\Curl \boldsymbol{E}=0$.
 \end{proof}
\end{lemma} 

Given a compatible strain field, the following lemma derives implications of the compatibility equation \eqref{compatcond} in terms of the restriction of the strain field on $\Omega - O$.
\begin{lemma}
\label{CompatibilityToLoop}
Given a simply connected region $\Omega\subset \mathbb{R}^2$ and a singular strain field $\boldsymbol{E}\in \mathcal{A}(\Omega,\Sym)$,
\begin{enumerate}[a.]
\item The compatibility condition \eqref{compatcond} implies
\begin{subequations}
\begin{align}
&\curl \curl (\boldsymbol{E}|_{\Omega-O})={0}~\text{and} \label{restrictedCompatibility_1} \\
&\int_{\partial B_\epsilon} \{ (\boldsymbol{E}|_{\Omega-O})(\boldsymbol{y})  + (\boldsymbol{y}-\boldsymbol{x})\times \curl ((\boldsymbol{E}|_{\Omega-O})(\boldsymbol{y})) \}\dy=\boldsymbol{0}\label{RestrictionCompatibility_2} 
\end{align}
\end{subequations}
for all $\boldsymbol{x} \in \Omega-O$.
\item For $\deg(\boldsymbol{E}) < -2$ the compatibility condition \eqref{compatcond} is equivalent to \eqref{restrictedCompatibility_1}.
\item  For $-2 \leq \deg(\boldsymbol{E}) < 0$ the compatibility condition \eqref{compatcond} is equivalent to \eqref{restrictedCompatibility_1} and \eqref{RestrictionCompatibility_2}.
\end{enumerate}
\begin{proof}
\textit{a.} This can be proved by following arguments from the proof of Lemma~\ref{strcompl3}.
\textit{b.} For $\boldsymbol{E}\in \mathcal{A}(\Omega,\Sym)$ with $\deg(\boldsymbol{E}) < -2$,  Lemma~\ref{ExistenceUniquenessExtensionLemma}\textit{a.} can be used to argue that $\Curl\Curl \boldsymbol{E}$ is a unique extension of $\curl \curl (\boldsymbol{E}|_{\Omega-O})$. Equation \eqref{restrictedCompatibility_1} would then immediately imply \eqref{compatcond} for the unique extension.
\textit{c.} According to Lemma~\ref{strcompl2} there exists an extension $\boldsymbol{E}$ of $\boldsymbol{E}|_{\Omega-O}$ such that $\Curl\Curl \boldsymbol{E}=0$. The required result is established by following Lemma~\ref{ExistenceUniquenessExtensionLemma}\textit{a.}  to argue the uniqueness of the extension given $\deg(\boldsymbol{E}) < 0$. 
\end{proof}
\end{lemma}

The preceding lemma is illustrated further by looking at specific examples. Consider an integrable strain field $\boldsymbol{E}$ such that $\Curl \Curl \boldsymbol{E}$ is also integrable, i.e., there exists integrable functions $\boldsymbol{f}:\Omega \to \Sym$ and ${g}:\Omega\to \mathbb{R}$ such that $\boldsymbol{E}(\boldsymbol{\phi})=\int_\Omega \langle \boldsymbol{f},\boldsymbol{\phi} \rangle \da$, for all $\boldsymbol{\phi}\in \mathcal{D}(\Omega,\Sym)$, and $\Curl \Curl \boldsymbol{E}( {\phi})=\int_\Omega  {g}{\phi}  \da$, for all ${\phi}\in \mathcal{D}(\Omega)$. Then,  $(\Curl\Curl \boldsymbol{E})|_{\Omega-O}={0}$ is equivalent to $\Curl\Curl \boldsymbol{E}={0}$. Also, if $\boldsymbol{E}\in \mathcal{A}(\Omega,\Sym)$ is of the above kind and $\boldsymbol{f}$ is smooth then $g=\curl \curl \boldsymbol{f}$. On the other hand, if we consider a strain field $\boldsymbol{E}\in \mathcal{A}(\Omega,\Sym)$, given by $\boldsymbol{E}(\boldsymbol{\phi})=\int_{\Omega}{\ln r}\left\langle  \boldsymbol{\phi} ,\boldsymbol{I}\right\rangle \da$, for arbitrary $\boldsymbol{\phi}\in \mathcal{D}(\Omega,\Sym)$, then $\curl\curl (\boldsymbol{E}|_{\Omega-O})={0}$ but $\Curl\Curl \boldsymbol{E} \neq {0}$. This is an example of a strain field $\boldsymbol{E}$ with $-2 \leq \deg(\boldsymbol{E}) < 0$, where although \eqref{restrictedCompatibility_1} is satisfied but \eqref{RestrictionCompatibility_2} is not, leading to a violation of the strain compatibility condition. If, however, we consider a strain field such that $\boldsymbol{E}(\boldsymbol{\phi})=\int_{\Omega}\left\langle  \boldsymbol{\phi} ,({\ln r} \boldsymbol{I} + \boldsymbol{e}_r\otimes \boldsymbol{e}_r) \right\rangle \da$, for which again $-2 \leq \deg(\boldsymbol{E}) < 0$, then both \eqref{restrictedCompatibility_1} and \eqref{RestrictionCompatibility_2} are satisfied, implying the compatibility of $\boldsymbol{E}$. In the last two examples, \eqref{restrictedCompatibility_1} along with \eqref{RestrictionCompatibility_2} are sufficient to establish the compatibility of the strain field. 
 
Consider a strain field $\boldsymbol{E}\in \mathcal{A}(\Omega,\Sym)$ such that $\deg(\boldsymbol{E})=\deg(\boldsymbol{E}|_{\Omega-O})$ and $\boldsymbol{E}|_{\Omega-O}$  behaves like $r^k$. Then, $\deg(\boldsymbol{E})=\deg(\boldsymbol{E}|_{\Omega-O})=-k-2$. For $k>0$, $\deg(\boldsymbol{E})<-2$, \eqref{restrictedCompatibility_1} is sufficient to establish the compatibility of strain fields. For $-2< k \leq 0$, $-2 \leq \deg(\boldsymbol{E}) < 0$, \eqref{restrictedCompatibility_1} and \eqref{RestrictionCompatibility_2} both should be satisfied for the strain compatibility. For $k\leq-2$, $\deg(\boldsymbol{E})\geq 0$, no conditions on $\boldsymbol{E}|_{\Omega-O}$ are sufficient to establish the strain compatibility, see Remark~\ref{stressoutinsuf}. In the following lemma, we provide the complete set of local conditions on $\boldsymbol{E}\in \mathcal{A}(\Omega,\Sym)$ which are equivalent to \eqref{compatcond}.
 \begin{lemma} \label{compstrcomplete}
For a simply connected region $\Omega\subset \mathbb{R}^2$ and a singular strain field $\boldsymbol{E}\in \mathcal{A}(\Omega,\Sym)$, the compatibility condition \eqref{compatcond} is equivalent to
\begin{subequations}
\begin{align}
&\curl \curl (\boldsymbol{E}|_{\Omega-O})={0}~\text{and} \label{restrictedCompatibility__1}\\
 &\boldsymbol{E}(\mathbb{A}\nabla^2 w^\alpha) =0 \label{restrictedCompatibility__2}
\end{align}
\end{subequations}
 for all $|\alpha|\leq\deg(\boldsymbol{E})+2$, where $\alpha$ is a multi-index, $\mathbb{A}$ is a linear map defined in Section~\ref{notation}, and $w^\alpha \in \mathcal{D}(\Omega)$ is as introduced in Lemma~\ref{RepresentationELemma}.
 
\begin{proof} Given \eqref{compatcond}, Equation~\eqref{restrictedCompatibility__1} follows by restricting it to $\Omega-O$ whereas Equations~\eqref{restrictedCompatibility__2} follows after using $w^\alpha$ as a test function in the distributional form of \eqref{compatcond}. On the other hand,
Equation \eqref{restrictedCompatibility__1} implies that $\Curl\Curl \boldsymbol{E}\in \mathcal{E}(\Omega,\mathbb{R}^2)$ with $\deg(\Curl\Curl \boldsymbol{E}) \leq \deg(\boldsymbol{E})+2$. Equation \eqref{restrictedCompatibility__2}, in conjunction with Lemma~\ref{RepresentationELemma}, then yields the desired result.
\end{proof} 
  \end{lemma}

Equation \eqref{restrictedCompatibility__2} is local at the singular point $O$ in the following sense. For any two singular strain fields $\boldsymbol{E}_1,\boldsymbol{E}_2\in \mathcal{A}(\Omega,\Sym)$ such that $\boldsymbol{E}_1|_\omega=\boldsymbol{E}_2|_\omega$, where $\omega \subset \Omega$ is an arbitrary open set such that $O\in \omega$, \eqref{restrictedCompatibility__2} is equivalent when written either in terms of $\boldsymbol{E}_1$ or $\boldsymbol{E}_2$. The compatibility of a singular strain field $\boldsymbol{E} \in \mathcal{A}(\Omega,\Sym)$, with $\deg(\boldsymbol{E})\geq -2$, is ensured if and only if both \eqref{restrictedCompatibility__1} and \eqref{restrictedCompatibility__2} are satisfied. Whenever $-2\leq\deg(\boldsymbol{E})<0$, \eqref{restrictedCompatibility__2} reduces to \eqref{RestrictionCompatibility_2}, and therefore compatibility can be checked through conditions only on  $\boldsymbol{E}|_{\Omega-O}$. Finally, only when $\deg(\boldsymbol{E})<-2$, \eqref{restrictedCompatibility__1} is sufficient for verifying the compatibility of the strain field.

 \subsection{Incompatibility of singular strain fields}
 
A strain field $\boldsymbol{E}\in \mathcal{D}'(\Omega,\Sym)$ is called incompatible if it does not satisfy the compatibility condition \eqref{compatcond}. Given a simply connected region $\Omega \subset \mathbb{R}^2$ and an incompatibility field $N\in \mathcal{D}'(\Omega)$, we say that a strain field $\boldsymbol{E}\in \mathcal{D}'(\Omega,\Sym)$ satisfies the incompatibility condition if 
\begin{equation}
\Curl \Curl \boldsymbol{E}=N. \label{incompatcond}
\end{equation}
We begin by considering point supported incompatibility fields. In the following lemma we obtain conditions on $N\in \mathcal{E}(\Omega)$ for there to exist a point supported strain field $\boldsymbol{E}  \in \mathcal{E}(\Omega,\Sym)$ which satisfies \eqref{incompatcond}. 
 
\begin{lemma}
\label{CurlCurlDipoleSource}
Consider $\Omega \subset \mathbb{R}^2$ and an incompatibility field $N \in \mathcal{E}(\Omega)$ with the representation 
\begin{equation}
N=\sum_{\alpha \in \mathbb{N}^2, |\alpha|\leq \deg(N)} N^\alpha \partial^\alpha \delta_O, ~\text{where}~N^\alpha \in \mathbb{R}. \label{representNE}
\end{equation}
If $ N^\alpha=0$ for $|\alpha|<2$ then there exists $\boldsymbol{E}\in \mathcal{E}(\Omega,\Sym)$ such that \eqref{incompatcond} holds true.
\begin{proof}
Given any multi-index $\alpha=(\alpha_1,\alpha_2)$, such that $|\alpha| \geq 2$, let $N_0= N^\alpha \partial^\alpha \delta_O.$ If $\alpha_1 \geq 2$ take $\boldsymbol{E}= N^\alpha \partial^{\alpha'} \delta_O \boldsymbol{e}_2\otimes\boldsymbol{e}_2$, where $\alpha'=(\alpha_1-2,\alpha_2)$. If $\alpha_1 = 0$ take $\boldsymbol{E}= N^\alpha \partial^{\alpha'} \delta_O \boldsymbol{e}_1\otimes \boldsymbol{e}_1$, where $\alpha'=(0,\alpha_2 -2)$.  If $\alpha_1 = 1$ take $\boldsymbol{E}= -(N^\alpha/2) \partial^{\alpha'} \delta_O (\boldsymbol{e}_2\otimes \boldsymbol{e}_1+\boldsymbol{e}_1\otimes \boldsymbol{e}_2)$, where $\alpha'=(0,\alpha_2 -1)$. For all the three cases we have $\Curl\Curl \boldsymbol{E}=N_0$. This establishes the lemma.
\end{proof}
\end{lemma}

In the next lemma we establish the implications on the restriction $\boldsymbol{E}|_{\Omega-O}$ given that $\boldsymbol{E}\in \mathcal{A}(\Omega,\Sym)$ satisfies \eqref{incompatcond} with a point supported incompatibility field $N\in \mathcal{E}(\Omega)$.
\begin{lemma}
\label{IncompatibilityCeasroIntegral}
Consider a simply connected region $\Omega \subset \mathbb{R}^2$, a singular strain field $\boldsymbol{E}\in \mathcal{A}(\Omega,\Sym)$, and an incompatibility field $N \in \mathcal{E}(\Omega)$ with the representation \eqref{representNE}. Then, 
\begin{enumerate}[a.]
\item The incompatibility equation  \eqref{incompatcond}
 implies
\begin{subequations} \label{RestrictionInCompatibility11}%
\begin{align}
&\curl \curl (\boldsymbol{E}|_{\Omega-O})={0}, \label{RestrictionInCompatibility_1} \\
&\int_{\partial B_\epsilon} \{ (\boldsymbol{E}|_{\Omega-O})(\boldsymbol{y})   + (\boldsymbol{y}-\boldsymbol{x})\times \curl ((\boldsymbol{E}|_{\Omega-O})(\boldsymbol{y})) \} \dy= -N^{(0,1)}\boldsymbol{e}_1+N^{(1,0)}\boldsymbol{e}_2 + N^{(0,0)} \boldsymbol{e}_3\times \boldsymbol{x},\label{RestrictionInCompatibility_2} \\
 &\deg(N) \leq \deg(\boldsymbol{E})+2, \label{RestrictionInCompatibility_0}
\end{align}
\end{subequations}
for all $\boldsymbol{x} \in \Omega-O$.
\item For $\deg(\boldsymbol{E}) < 0$, the incompatibility condition \eqref{incompatcond} is equivalent to Equations \eqref{RestrictionInCompatibility11}.
\end{enumerate}
\begin{proof}
\textit{a.} Equation~\eqref{incompatcond}, with $N \in \mathcal{E}(\Omega)$, restricted to $\Omega -O$ gives  \eqref{RestrictionInCompatibility_1}. Equation~\eqref{RestrictionInCompatibility_0} follows from point~v. given at the end of Section~\ref{sddeg}. Let $\boldsymbol{E}_2 =(1/2\pi r) \left( -N^{(0,1)}\boldsymbol{e}_1+N^{(1,0)}\boldsymbol{e}_2 \right)\otimes \boldsymbol{e}_{\theta}+ N^{(0,0)}r\ln r \boldsymbol{I}$. Then, according to Lemma~\ref{CurlCurlDipoleSource}, there exists $\boldsymbol{E}_1\in \mathcal{E}(\Omega,\Sym)$ such that the strain field $\boldsymbol{E}_3=\boldsymbol{E}_1+\boldsymbol{E}_2$ satisfies $\Curl \Curl \boldsymbol{E}_3=N$. Consequently $\Curl \Curl (\boldsymbol{E}-\boldsymbol{E}_3)=0$ which, using Lemma~\ref{CompatibilityToLoop}, leads to \eqref{RestrictionInCompatibility_2}.
\textit{b.} Equation~\eqref{RestrictionInCompatibility_0} with $\deg(\boldsymbol{E}) < 0$ implies that $\deg(N)<2$. This along with Lemma~\ref{RepresentationELemma} and Equations~\eqref{RestrictionInCompatibility_1} and \eqref{RestrictionInCompatibility_2} proves \eqref{incompatcond}.
\end{proof}
\end{lemma}

Considering the form of Equation~\eqref{RestrictionInCompatibility_2}, we can interpret $-N^{(0,1)}\boldsymbol{e_1}+N^{(1,0)}\boldsymbol{e_2}$ as the Burgers vector of an isolated dislocation and $N^{(0,0)}$ as the disclination charge of an isolated disclination at point $O$. The incompatibility field $N\in \mathcal{E}(\Omega)$ with $\deg(N)<2$ therefore captures the point supported incompatibility due to an isolated dislocation and disclination placed at point $O$ \cite{van2012distributional}. The point supported incompatibility fields $N\in \mathcal{E}(\Omega)$ with higher degrees of divergence imply the presence of other point defects such as a dislocation dipole (or quadrupole, etc.), extra-matter, vacancy, or a concentrated heat/growth source.

The incompatibility of a strain field $\boldsymbol{E}\in \mathcal{A}(\Omega,\Sym)$ can be completely characterized in terms of incompatibility of the restricted field $\boldsymbol{E}|_{\Omega-O}$ only for incompatibility fields with $\deg(N)<2$. As discussed earlier in the context of stress and compatible strain, see Remark~\ref{stressoutinsuf}, the restricted field will in general be never sufficient to describe the incompatibility of the strain field. In the following lemma, we derive the necessary and sufficient local conditions associated with a singular strain field which are equivalent to \eqref{incompatcond}.

 \begin{lemma} \label{incompstrcomplete}
Consider a simply connected region $\Omega \subset \mathbb{R}^2$, a singular strain field $\boldsymbol{E}\in \mathcal{A}(\Omega,\Sym)$, and an incompatibility field $N \in \mathcal{A}(\Omega)$. The incompatibility equation  \eqref{incompatcond} is equivalent to
\begin{subequations}  \label{RestrictionInCompatibility22}%
\begin{align}
&\curl \curl (\boldsymbol{E}|_{\Omega-O})={N}|_{\Omega-O},\label{restrictedInCompatibility__1} \\
& \boldsymbol{E}(\mathbb{A}\nabla^2 w^\alpha) =N(w^\alpha),~\text{and} \label{restrictedInCompatibility__2} \\
&\deg(N) \leq \deg (\boldsymbol{E})+2,  \label{restrictedInCompatibility__3}
\end{align}
\end{subequations}
 for all $|\alpha|\leq\deg(\boldsymbol{E})+2$, where $\alpha$ is a multi-index, $\mathbb{A}$ is a linear map defined in Section~\ref{notation}, and $w^\alpha \in \mathcal{D}(\Omega)$ is as introduced in Lemma~\ref{RepresentationELemma}.
 \begin{proof}
Given \eqref{incompatcond}, Equation~\eqref{restrictedInCompatibility__1} follows by restricting it to $\Omega-O$, \eqref{restrictedInCompatibility__3} follows from point~v. given at the end of Section~\ref{sddeg}, and \eqref{restrictedInCompatibility__2} is obtained on using $w^\alpha$ as a test function in the distributional form of \eqref{incompatcond}. On the other hand, 
Equation \eqref{restrictedInCompatibility__1} implies that $(\Curl\Curl \boldsymbol{E}-N)\in \mathcal{E}(\Omega)$. Due to \eqref{restrictedInCompatibility__3}, $\deg(\Curl\Curl \boldsymbol{E}-N) \leq \deg(\boldsymbol{E})+2$. Equation~\eqref{restrictedInCompatibility__2}, in conjunction with Lemma~\ref{RepresentationELemma}, subsequently yields \eqref{incompatcond}.
\end{proof} 
 \end{lemma}
 
 Equation~\eqref{restrictedInCompatibility__2} is local at point $O$ in the same sense as elaborated at the end of Section~\ref{compsingstr}. It is clear that for a strain field $\boldsymbol{E}\in \mathcal{A}(\Omega,\Sym)$, such that $\deg(\boldsymbol{E})<-2$, and an incompatibility field $N \in \mathcal{A}(\Omega)$, such that $\deg(N)<0$,  \eqref{restrictedInCompatibility__1} is equivalent to the incompatibility equation  \eqref{incompatcond}. For a strain field with 
$\deg(\boldsymbol{E})\geq-2$, additional conditions of the form of \eqref{restrictedInCompatibility__2} are required for establishing the equivalence.  In particular, if $-2\leq \deg(\boldsymbol{E})<0$ and $N$ is a point supported incompatibility field, such that $\deg(N)<2$, \eqref{restrictedInCompatibility__2}  reduces to \eqref{RestrictionInCompatibility_2}. Note that both \eqref{restrictedInCompatibility__1} and \eqref{RestrictionInCompatibility_2} are conditions 
on $\boldsymbol{E}|_{\Omega-O}$. For strain fields with non-negative degree of divergence, however, the incompatibility condition \eqref{incompatcond} cannot be interpreted only in terms of $\boldsymbol{E}|_{\Omega-O}$. Indeed, consider a strain field $\boldsymbol{E}$ with $\deg(\boldsymbol{E}) \geq 0$. The field $\boldsymbol{E}_1$, defined as $\boldsymbol{E}_1=\boldsymbol{E}+\delta_O \boldsymbol{I}$, satisfies
${\boldsymbol{E}_1}|_{\Omega-O}=\boldsymbol{E}|_{\Omega-O}$ but $\Curl \Curl \boldsymbol{E}_1 \neq \Curl \Curl \boldsymbol{E}$.

\begin{remark} (Sources of strain incompatibility)
Dislocations, disclinations, and non-uniform temperature fields are possible sources of strain incompatibility \cite{de1981view}. For instance, given a dislocation density field $\boldsymbol{A}\in \mathcal{D}'(\Omega,\mathbb{R}^2)$ and a disclination density field ${\Theta}\in \mathcal{D}'(\Omega)$, the incompatibility field $N\in \mathcal{D}'(\Omega)$ can be written as 
$N= \Curl \boldsymbol{A} + {\Theta}$ \cite{pandey2020topological}.
In the absence of defects but given a temperature field $\vartheta \in \mathcal{D}'(\Omega)$, the induced incompatibility field  is $N=\Delta \vartheta$. In particular, if $\vartheta=\ln r$ then $N=\delta_O$. If $\boldsymbol{A}\in \mathcal{A}(\Omega,\mathbb{R}^2)$, ${\Theta}\in \mathcal{A}(\Omega)$, and $\vartheta \in \mathcal{A}(\Omega)$ then $N \in \mathcal{A}(\Omega)$. On the other hand, if $\boldsymbol{A}\in \mathcal{E}(\Omega,\mathbb{R}^2)$, ${\Theta}\in \mathcal{E}(\Omega)$, and $\vartheta \in \mathcal{E}(\Omega)$ then $N \in \mathcal{E}(\Omega)$. Given an incompatibility field $N$ coming from dislocations and disclinations, we do not have a unique prescription for $\boldsymbol{A}$ or ${\Theta}$.
For example, incompatibility $N=\delta_O$ can result from either of the two pairs: $\boldsymbol{A}_1=\boldsymbol{0}$, ${\Theta_1}=\delta_O$ and $\boldsymbol{A}_2(\boldsymbol{\phi})=\int_\Omega({1}/{2\pi r})\left\langle \boldsymbol{\phi}, \boldsymbol{e}_\theta \right\rangle \da$, for all $\boldsymbol{\phi}\in \mathcal{D}(\Omega,\mathbb{R}^2)$, ${\Theta_2}={0}$.
In the case of vanishing disclination density, i.e., $\Theta =0$, there exists a distortion field $\boldsymbol{\beta}\in \mathcal{D}'(\Omega,\Lin)$ such that $\Curl \boldsymbol{\beta}=\boldsymbol{A}$ and $\Curl \Curl \boldsymbol{E}=N$, where $\boldsymbol{E}=(1/2)(\boldsymbol{\beta}+\boldsymbol{\beta}^T )$ \cite{pandey2020topological}. It is clear that the existence of $\boldsymbol{\beta}$ is undetermined by solely prescribing an incompatibility field $N$. \label{incomdefect}
\end{remark}

\section{The stress problem of linear elasticity} \label{uniqueness}
\subsection{Uniqueness of the stress solution}

The stress problem of linear incompatible elasticity is a traction boundary value problem for the determination of stress field in the elastic body for a given distribution of forces (body forces and traction) and incompatibility. In this section we will state the boundary value problem in a distributional form and prove that the resulting solution is unique. The uniqueness result significantly generalizes the earlier work by Sternberg and coauthors \cite{sternberg1955concept, Turteltaub1968OnCL, gurtin1973linear}, as elaborated towards the end of the section. Consider a simply connected region $\Omega \subset \mathbb{R}^2$ with a given prescription of a singular body force field $\boldsymbol{B}\in\mathcal{A}(\Omega,\mathbb{R}^2)$, a smooth traction field $\boldsymbol{t}:\partial \Omega \to \mathbb{R}^2$ on the boundary $\partial \Omega$ of $\Omega$ (with outward normal $\boldsymbol{n}$), and a singular incompatibility field ${N}\in\mathcal{A}(\Omega)$. The stress field $\boldsymbol{\sigma} \in \mathcal{A}(\Omega,\Sym)$ is determined by solving the following boundary value problem:
\begin{subequations} \label{tracbvp}
\begin{align}
&\Div \boldsymbol{\sigma}+\boldsymbol{B} = \boldsymbol{0},\label{tracbvp1} \\
&\Curl \Curl \mathbb{C}^{-1}\boldsymbol{\sigma} = {N},~\text{and} \label{tracbvp2}\\
 &\boldsymbol{\sigma}\boldsymbol{n} = \boldsymbol{t}~\text{on}~\partial\Omega,  \label{tracbvp3}
\end{align}
\end{subequations}
where $\mathbb{C}:\Sym \to \Sym$ is the elasticity tensor such that $\boldsymbol{\sigma}=\mathbb{C}\boldsymbol{E}$ with $\boldsymbol{E}\in \mathcal{A}(\Omega,\Sym)$ representing the elastic strain field. We assume $\mathbb{C}$ to be symmetric and positive-definite, i.e., $\langle \mathbb{C}\boldsymbol{V}_1, \boldsymbol{V}_2\rangle= \langle \mathbb{C}\boldsymbol{V}_2, \boldsymbol{V}_1\rangle$, for all $\boldsymbol{V}_1,\boldsymbol{V}_2 \in \Sym$, and $\langle \mathbb{C}\boldsymbol{V}, \boldsymbol{V}\rangle > 0$, for all $\boldsymbol{V}\in \Sym$ such that $\boldsymbol{V} \neq \boldsymbol{0}$. Note that the boundary condition \eqref{tracbvp3} is a pointwise condition since the restriction of $\boldsymbol{\sigma}$ to $\Omega - O$ is smooth and $O$ lies in the interior of the domain. The following lemma establishes the uniqueness of a stress solution to the boundary value problem. A unique stress field would immediately imply a unique elastic strain field. However, the existence of a unique displacement solution $\boldsymbol{u}\in \mathcal{A}(\Omega)$, modulo translation and rotation, such that $\boldsymbol{E} = (1/2)(\nabla \boldsymbol{u}+(\nabla \boldsymbol{u})^T)$, can be argued only when the incompatibility field $N$ is identically zero.

\begin{lemma}
\label{LinearElasticityUniquenessLemma}%
For a simply connected region $\Omega\subset \mathbb{R}^2$ and a symmetric,  positive-definite elasticity tensor $\mathbb{C}$, let $\boldsymbol{\sigma}_1,\boldsymbol{\sigma}_2 \in \mathcal{A}(\Omega,\Sym)$ be any two solutions of the problem \eqref{tracbvp}, with a given distribution of body force $\boldsymbol{B}\in\mathcal{A}(\Omega,\mathbb{R}^2)$, smooth traction $\boldsymbol{t}:\partial \Omega \to \mathbb{R}^2$, and incompatibility ${N}\in\mathcal{A}(\Omega)$. Then, $\boldsymbol{\sigma}_1=\boldsymbol{\sigma}_2$.
\begin{proof}
The field $\bar{\boldsymbol{\sigma}} =  \boldsymbol{\sigma}_1-\boldsymbol{\sigma}_2$ satisfies
\begin{subequations} \label{tracbvphom}
\begin{align}
&\Div \bar{\boldsymbol{\sigma}} = \boldsymbol{0},\label{tracbvphom1} \\
&\Curl \Curl \mathbb{C}^{-1}\bar{\boldsymbol{\sigma}} = 0,~\text{and} \label{tracbvphom2}\\
 &\bar{\boldsymbol{\sigma}} \boldsymbol{n} = \boldsymbol{0}~\text{on}~\partial\Omega.  \label{tracbvphom3}
\end{align}
\end{subequations}
Equation~\eqref{tracbvphom1} implies that there exist a scalar distribution $\Phi \in \mathcal{A}(\Omega)$ such that $\bar{\boldsymbol{\sigma}} = \mathbb{A}\nabla^2 \Phi$, where $\mathbb{A}$ is a linear map defined in Section~\ref{notation}. Substituting $\bar{\boldsymbol{\sigma}}$ in terms of $\Phi$ into Equation~\eqref{tracbvphom2} yields
\begin{equation}
\label{EquationAiryStressUniquenessLemma}
\Curl\Curl (\mathbb{C}^{-1}\mathbb{A}\nabla^2 \Phi ) =0.
\end{equation}
We note the following two identities: (i) 
 For $\boldsymbol{T}\in \mathcal{D}'(\Omega,\Sym)$,  $\Curl \Curl \boldsymbol{T} = \Div\Div (\mathbb{A} \boldsymbol{T})$ and (ii) For $T \in \mathcal{D}'(\Omega)$, and a fixed $\mathbb{D} \in (\mathbb{R}^2)^4$, $\Div\Div (\mathbb{D} \nabla^2 {T})= \left\langle \mathbb{D}, \nabla^4 T  \right\rangle$. Using these, Equation~\eqref{EquationAiryStressUniquenessLemma} can be rewritten as
\begin{equation}
 \left\langle \mathbb{A} \mathbb{C}^{-1}\mathbb{A} , \nabla^4 \Phi  \right\rangle = 0.
\end{equation}
Using the identity
\begin{equation}
\left\langle \mathbb{A} \mathbb{C}^{-1}\mathbb{A} ,\boldsymbol{v}\otimes \boldsymbol{v}\otimes \boldsymbol{v}\otimes \boldsymbol{v}  \right\rangle = \left\langle  \mathbb{C}^{-1}  (\boldsymbol{e}_3 \times\boldsymbol{v})\otimes (\boldsymbol{e}_3 \times\boldsymbol{v}),(\boldsymbol{e}_3 \times\boldsymbol{v})\otimes (\boldsymbol{e}_3 \times\boldsymbol{v})  \right\rangle,
\end{equation}
for any $\boldsymbol{v} \in \mathbb{R}^2$, we assert that the ellipticity of the operator $\Curl\Curl (\mathbb{C}^{-1}\mathbb{A}\nabla\nabla)$ follows from the positive-definiteness of $\mathbb{C}$. Lemma~\ref{EllipticRegularityLemma} can then be used to conclude that $\singsupp(\Phi)=\singsupp(0)$. Hence $\Phi$, and therefore $\bar{\boldsymbol{\sigma}}$, is a smooth field. The smoothness of $\bar{\boldsymbol{\sigma}}$, in conjunction with the classical uniqueness theorem for the smooth stress problem in linear elasticity \cite{gurtin1973linear}, implies $\bar{\boldsymbol{\sigma}}=\boldsymbol{0}$.
\end{proof}

\end{lemma}

The uniqueness of the stress solution is a crucial property in the linear elasticity framework. The existence of non-unique solutions would yield non-trivial stress fields as solutions in response to vanishing source fields. Lemma~\ref{LinearElasticityUniquenessLemma} establishes the uniqueness of a solution for the traction  boundary value problem with stress fields singular at an isolated point in $\Omega$. The uniqueness theorem extends naturally to more general singular stress fields, for example to fields with multiple points of singularity.

In their work Sternberg and coauthors \cite{sternberg1955concept, Turteltaub1968OnCL, gurtin1973linear} have formulated a stress boundary value problem with the singular stress satisfying certain scaling assumptions in the vicinity of the singular points. The scaling assumptions are in fact such that the degree of divergence of the stress field, with respect to the singular point, remains negative. In the proposed formalism, the stress field is a smooth map away from the points of singularity and the body force field is such that it allows for concentrated loads at the singular points. There is no consideration of the incompatibility field. A uniqueness theorem for the singular stress solution, satisfying the scaling assumptions, for the considered boundary value problem is also proved. More importantly, it is mentioned that the solutions become non-unique when the stress fields do not satisfy the scaling assumptions \cite{sternberg1955concept}. This is illustrated in detail for the stress solution in response to a force dipole. The reason for this apparent lack of unique solutions in their framework is due to the limited consideration of stress field only as a map away from the singular point. In other words, considerations are limited to $\boldsymbol{\sigma}|_{\Omega-O}$ rather than the full stress field $\boldsymbol{\sigma}$. Under the stronger scaling assumptions, when $\deg(\boldsymbol{\sigma}|_{\Omega-O}) <0$ (which implies a unique extension $\boldsymbol{\sigma}$ such that $\deg(\boldsymbol{\sigma}) < 0$, see Lemma~\ref{ExistenceUniquenessExtensionLemma}\textit{a.}), we have established in Lemma~\ref{stressequiv2} that stress equilibrium condition can be equivalently written in terms of conditions on the restricted map $\boldsymbol{\sigma}|_{\Omega-O}$. The uniqueness in such cases follows from considerations of $\boldsymbol{\sigma}|_{\Omega-O}$ alone. On the other hand, this is not so whenever $\deg(\boldsymbol{\sigma}) \geq 0$, as is in the case of a force dipole, and one is required to consider conditions on the full stress field $\boldsymbol{\sigma}$. As stated in Lemma~\ref{LinearElasticityUniquenessLemma}, the stress solutions to the boundary value problem \eqref{tracbvp} of incompatible linear elasticity always satisfy the property of uniqueness when the stress field is considered as a singular distribution on $\Omega$ without any a priori scaling assumptions on the stress field. The examples in the following section further expand our point of view.

\subsection{A general stress solution}

Let the elasticity tensor $\mathbb{C}$  be given as for an isotropic plane strain scenario, i.e., $\boldsymbol{E} = ((1+\nu)/E) \boldsymbol{\sigma} - (\nu(1+\nu)/E) (\tr \boldsymbol{\sigma}) \boldsymbol{I}$, where $E$ is the Young's modulus and $\nu$ is the Poisson's ratio.  The stress field
\begin{equation}
\boldsymbol{\sigma}_1(\boldsymbol{\phi})=\frac{1-2\nu}{4\pi(1-\nu)}\int_\Omega \left\langle \frac{1}{r}\left(\frac{2\nu-3}{1-2\nu} \cos \theta \boldsymbol{e}_r\otimes\boldsymbol{e}_r+ \sin \theta (\boldsymbol{e}_r\otimes\boldsymbol{e}_\theta+\boldsymbol{e}_\theta\otimes\boldsymbol{e}_r)+ \cos \theta \boldsymbol{e}_\theta \otimes\boldsymbol{e}_\theta \right), \boldsymbol{\phi}\right\rangle \da, \nonumber
\end{equation}
  for all $\boldsymbol{\phi}\in \mathcal{D}(\Omega,\Lin)$, is a solution to the following pair of equations:
  \begin{equation}
  \Div \boldsymbol{\sigma}_1+\delta_O\boldsymbol{e}_1=\boldsymbol{0}~\text{and}~
  \Curl \Curl \mathbb{C}^{-1}\boldsymbol{\sigma}_1=0,
  \end{equation}
  whereas the stress field 
\begin{equation}
 \boldsymbol{\sigma}_2(\boldsymbol{\phi})=\frac{E}{8\pi(1-\nu^2)}\int_\Omega \left\langle \left( (2\ln r +1) \boldsymbol{e}_r\otimes\boldsymbol{e}_r+ (2\ln r +3) \boldsymbol{e}_\theta \otimes\boldsymbol{e}_\theta \right), \boldsymbol{\phi}\right\rangle \da,
 \end{equation}
for all $\boldsymbol{\phi}\in \mathcal{D}(\Omega,\Lin)$, is a solution to the following pair of equations:
  \begin{equation}
  \Div \boldsymbol{\sigma}_2=\boldsymbol{0}~\text{and}~
  \Curl \Curl \mathbb{C}^{-1} \boldsymbol{\sigma}_2=\delta_O.
  \end{equation}
Both of these stress field solutions are integrable with $\deg(\boldsymbol{\sigma})<0$. The restricted map $\boldsymbol{\sigma}|_{\Omega-O}$, away from the point of singularity,  is sufficient to describe the complete solution and its uniqueness for the corresponding boundary value problems. For weaker point supported body force and incompatibility fields, i.e. those satisfying $\deg(\boldsymbol{B})\geq 1$ or $N  \geq 2$, we have $\deg(\boldsymbol{\sigma})\geq 0$. In such cases the restricted map no longer  determines the complete solution. The two stress fields $\boldsymbol{\sigma}_1$ and $\boldsymbol{\sigma}_2$ can be used to construct more general solutions. 
For example, given point supported body force and incompatibility fields such that 
\begin{equation}
\boldsymbol{B}=\sum_{\alpha \in \mathbb{N}^2, |\alpha|\leq\deg(\boldsymbol{B})} B^\alpha \partial^\alpha \delta_O \boldsymbol{e}_1~\text{and}~N=\sum_{\alpha \in \mathbb{N}^2, |\alpha|\leq\deg(N)} N^\alpha \partial^\alpha \delta_O,
\end{equation}
 the stress field 
\begin{equation}
\label{LinearElasticitySolutionPointSource}
\boldsymbol{\sigma}=\sum_{\alpha \in \mathbb{N}^2, |\alpha|<\deg(\boldsymbol{B})} B^\alpha \partial^\alpha \boldsymbol{\sigma}_1 +  \sum_{\alpha \in \mathbb{N}^2, |\alpha|<\deg(N)} N^\alpha \partial^\alpha \boldsymbol{\sigma}_2
\end{equation} 
satisfies
\begin{equation}
  \Div \boldsymbol{\sigma} + \boldsymbol{B}=\boldsymbol{0}~\text{and}~
  \Curl \Curl \mathbb{C}^{-1} \boldsymbol{\sigma}=N. \label{stresscompgovernex}
  \end{equation} 
The stress solution for a complete boundary value problem can be obtained by superposing this solution with the smooth stress solutions corresponding to trivial bulk sources and appropriate traction fields. We emphasize that in Equation~\eqref{LinearElasticitySolutionPointSource} the derivatives have to be necessarily interpreted as distributional derivatives in order to obtain the complete solution. Only the restricted field away from the singularity can be obtained using the smooth derivative in $\Omega-O$.

Given $\boldsymbol{B}=\boldsymbol{0}$ and $N=(2(\nu^2-1)/E) \Delta \delta_O$, the stress field
\begin{equation}
\boldsymbol{\sigma}_3(\boldsymbol{\phi})=\lim_{\epsilon \to 0} \int_{\Omega-B_\epsilon} \left\langle \frac{1}{\pi r^2}\left(-\boldsymbol{e}_r \otimes\boldsymbol{e}_r+\boldsymbol{e}_\theta \otimes\boldsymbol{e}_\theta \right), \boldsymbol{\phi}\right\rangle \da-\langle \boldsymbol{\phi}(O),\boldsymbol{I} \rangle,
\end{equation}
for all $\boldsymbol{\phi} \in \mathcal{D}(\Omega,\Lin)$, satisfies \eqref{stresscompgovernex}. On the other hand, when $\boldsymbol{B}=(2(\nu-1)/(1-2\nu))\nabla \delta_O$ and $N=0$, the stress field
\begin{equation}
\boldsymbol{\sigma}_4(\boldsymbol{\phi})=\lim_{\epsilon \to 0} \int_{\Omega-B_\epsilon} \left\langle \frac{1}{\pi r^2}\left(-\boldsymbol{e}_r\otimes\boldsymbol{e}_r+\boldsymbol{e}_\theta\otimes\boldsymbol{e}_\theta\right), \boldsymbol{\phi}\right\rangle \da+\frac{1}{1-2\nu}\langle \boldsymbol{\phi}(O),\boldsymbol{I} \rangle,
\end{equation}
for all $\boldsymbol{\phi} \in \mathcal{D}(\Omega,\Lin)$, satisfies \eqref{stresscompgovernex}. The stress fields $\boldsymbol{\sigma}_3$ and $\boldsymbol{\sigma}_4$ are such that their restrictions to $\Omega - O$ are identical but they are different solutions in response to different source fields. For both the solutions, $\deg(\boldsymbol{\sigma}|_{\Omega-O})=0$. The restricted field can therefore have non-unique extensions in $\Omega$, as established in Lemma~\ref{ExistenceUniquenessExtensionLemma}\textit{b.} This is an instructive example which clearly demonstrates the need for considering the full stress solution rather than its restriction over the domain outside the point of singularity.

\section{Force on a defect} \label{forceonadefect}

For the purpose of this section we assume the incompatibility field, over a simply connected region $\Omega \subset \mathbb{R}^2$, to be given solely in terms of dislocation density $\boldsymbol{A} \in  \mathcal{A}(\Omega,\mathbb{R}^2)$, i.e., $N= \Curl \boldsymbol{A}$, see Remark~\ref{incomdefect}. The disclination density $\Theta$ is assumed to vanish throughout.
Given $\boldsymbol{A}$, the induced incompatibility $N\in  \mathcal{A}(\Omega)$ can be used in the boundary value problem \eqref{tracbvp} to obtain a unique elastic strain field $\boldsymbol{E}\in  \mathcal{A}(\Omega,\Sym)$. Noting  the relation between the strain field and the dislocation density, $\boldsymbol{A}=\Curl \boldsymbol{E}+\boldsymbol{e}_3\otimes \nabla a$ \cite{pandey2020topological}, where $a \in  \mathcal{A}(\Omega)$, we can deduce the following form for the elastic distortion field $\boldsymbol{\beta}\in  \mathcal{A}(\Omega,\Lin)$:
\begin{equation}
\boldsymbol{\beta}=\boldsymbol{E} + a(\boldsymbol{e}_1 \otimes \boldsymbol{e}_2 -\boldsymbol{e}_2 \otimes \boldsymbol{e}_1).
\end{equation}
The elastic stress-strain relation $\boldsymbol{\sigma} = \mathbb{C}\boldsymbol{E}$ can be equivalently written as $\boldsymbol{\sigma} = \mathbb{C}\boldsymbol{\beta}$ by extending $\mathbb{C}: \Sym \to \Sym$ to a linear map $\mathbb{C}: \Lin \to \Sym$ such that $\mathbb{C} \boldsymbol{V} = \boldsymbol{0}$, for any $\boldsymbol{V}\in \Skw$.

The notion of a force acting on a defect, introduced originally by Eshelby \cite{eshelby1951force, eshelby1956continuum}, is essentially a thermodynamic concept which is related to the change in the total free energy of the elastic domain as the defect moves by an infinitesimal distance in the domain. It plays a central role in dealing with the problems of defect equilibrium and defect kinetics. Our aim is to demonstrate the utility of the methods developed in the present work for rigorously deriving the expressions for the force acting on isolated defects in flat two-dimensional domains. Our derivations in fact give a generalized force which is related to both translational and non-translational changes in the defect configuration. Toward this end, let distortion $\boldsymbol{\beta}_1 \in  \mathcal{A}(\Omega,\Lin)$ and stress $\boldsymbol{\sigma}_1\in  \mathcal{A}(\Omega,\Sym)$ be the solution to a boundary value problem with $\boldsymbol{B} = \boldsymbol{0}$, $\boldsymbol{t}=\boldsymbol{0}$, and an incompatibility field given in terms of a known dislocation density $\boldsymbol{A} \in  \mathcal{A}(\Omega,\mathbb{R}^2)$. Let distortion $\boldsymbol{\beta}_2 \in C^\infty (\Omega,\Lin)$ and stress $\boldsymbol{\sigma}_2 \in C^\infty (\Omega,\Sym)$ be the smooth solution to another boundary value problem with $\boldsymbol{B} = \boldsymbol{0}$, $N={0}$, but a known smooth traction field $\boldsymbol{t}$ on $\partial \Omega$. The solution to the latter boundary value problem can also be understood as a response due to sources which are external to the domain. We define the interaction Eshelby tensor field $\boldsymbol{J}^{I} \in \mathcal{A}(\Omega,\Lin)$ as
\begin{equation}
\boldsymbol{J}^{I}=\frac{1}{2}\left\langle\mathbb{C}\boldsymbol{\beta}_1,\boldsymbol{\beta}_2\right\rangle\boldsymbol{I}-{\boldsymbol{\beta}_2}^T\boldsymbol{\sigma}_1+\frac{1}{2}\left\langle\mathbb{C}\boldsymbol{\beta}_2,\boldsymbol{\beta}_1\right\rangle\boldsymbol{I}-{\boldsymbol{\beta}_1}^T\boldsymbol{\sigma}_2. \label{eshtenint}
\end{equation}
We introduce $\boldsymbol{F}^I = \Div \boldsymbol{J}^{I}$ as the generalized force on the defect due to the interaction between the singular field $\boldsymbol{\beta}_1$ and the smooth field $\boldsymbol{\beta}_2$. With considerations made above, the distortion and stress fields satisfy $\Curl \boldsymbol{\beta}_1=\boldsymbol{A}$, $\Div \boldsymbol{\sigma}_1=\boldsymbol{0}$, $\Curl \boldsymbol{\beta}_2=\boldsymbol{0}$, and $\Div \boldsymbol{\sigma}_2=\boldsymbol{0}$, yielding
\begin{equation}
\label{GeneralizedDrivingForceonDistribution}
\boldsymbol{F}^I = (\boldsymbol{\sigma}_2\boldsymbol{A})\times \boldsymbol{e}_3.
\end{equation}  
The definition of the interaction tensor and the generalized force are both motivated from Eshelby's work \cite{eshelby1951force, eshelby1956continuum}. To elaborate we assume, for now, that $\boldsymbol{\beta}_1$ and $\boldsymbol{\sigma}_1$ are smooth fields. Let $\boldsymbol{\beta} = \boldsymbol{\beta}_1 + \boldsymbol{\beta}_2$; the corresponding stress $\boldsymbol{\sigma}=\mathbb{C}\boldsymbol{\beta}$ satisfies $\boldsymbol{\sigma} = \boldsymbol{\sigma}_1 + \boldsymbol{\sigma}_2$. The Eshelby tensor is given by $\boldsymbol{J}=W\boldsymbol{I}-\boldsymbol{\beta}^T\boldsymbol{\sigma}$, where $W=({1}/{2})\left\langle\mathbb{C}\boldsymbol{\beta},\boldsymbol{\beta}\right\rangle $ is the strain energy density. It can be decomposed as $\boldsymbol{J}=\boldsymbol{J}_1+\boldsymbol{J}_2+\boldsymbol{J}^{I}$,
where $\boldsymbol{J}_1=\frac{1}{2}\left\langle\mathbb{C}\boldsymbol{\beta}_1,\boldsymbol{\beta}_1\right\rangle\boldsymbol{I}-{\boldsymbol{\beta}_1}^T\boldsymbol{\sigma}_1$ and $\boldsymbol{J}_2=\frac{1}{2}\left\langle\mathbb{C}\boldsymbol{\beta}_2,\boldsymbol{\beta}_2\right\rangle\boldsymbol{I}-{\boldsymbol{\beta}_2}^T\boldsymbol{\sigma}_2$ are self fields associated with $\boldsymbol{\beta}_1$ and $\boldsymbol{\beta}_2$, respectively. The term $\boldsymbol{J}^{I}$, which is of the form \eqref{eshtenint}, represents the interaction between the two fields.

We note two technical points before moving on to evaluate forces on some specific defect configurations. The definitions of the interaction Eshelby tensor field and the generalized force on the defect are local in nature and so are the regularity assumptions on distortion fields $\boldsymbol{\beta}_1$ and $\boldsymbol{\beta}_2$. This allows us to consider $\boldsymbol{\beta}_2$ fields which are also singular, as long as $\singsupp(\boldsymbol{\beta}_1) \cap \singsupp(\boldsymbol{\beta}_2)=\emptyset$. There will always exist a neighbourhood $\omega$ of $\singsupp(\boldsymbol{\beta}_1)$ such that $\boldsymbol{\beta}_2|_{\omega}$ is smooth and both $\boldsymbol{J}^I$ and $\boldsymbol{F}^I$ are well defined on $\omega$. Secondly, 
consider a sequence of distortion fields, $\boldsymbol{\beta}_1^j \in \mathcal{A}(\Omega,\Lin)$ such that $\boldsymbol{\beta}_1^j\to \boldsymbol{\beta}_1$. Let ${\boldsymbol{J}^I}^j \in \mathcal{A}(\Omega,\Lin)$ be the interaction Eshelby tensor associated with the interaction between $\boldsymbol{\beta}_1^j$ and $\boldsymbol{\beta}_2$ fields. Then, ${\boldsymbol{J}^I}^j\to \boldsymbol{J}^I$ and $\Div{\boldsymbol{J}^I}^j\to \Div \boldsymbol{J}$, as a consequence of the linear dependence of ${\boldsymbol{J}^I}^j$ on $\boldsymbol{\beta}_1^j$. 

We now provide several examples for point supported dislocation density fields, i.e., $\boldsymbol{A}\in \mathcal{E}(\Omega,\mathbb{R}^2)$. The generalized force too is then point supported, i.e., $\boldsymbol{F}^I \in \mathcal{E}(\Omega,\mathbb{R}^2)$. Since $\deg(\boldsymbol{F}^I) \leq \deg(\boldsymbol{A})$, we use Lemma~\ref{RepresentationELemma} to write the following representation of the generalized force:
\begin{equation}
{\boldsymbol{F}^I}=\sum_{\alpha \in \mathbb{N}^2, |\alpha|\leq\deg(\boldsymbol{A})} {\boldsymbol{F}^I}^\alpha \partial^\alpha \delta_O, \label{represdrivforce}
\end{equation}
where ${\boldsymbol{F}^I}^\alpha \in \mathbb{R}^2$. In particular, since $\boldsymbol{F}^I=\Div \boldsymbol{J}^I$ and ${\boldsymbol{F}^I}|_{\Omega-O}=\boldsymbol{0}$, we identify
\begin{equation}
{\boldsymbol{F}^I}^{(0,0)}=\int_{\partial B_\epsilon} \boldsymbol{J}^I \boldsymbol{n} \dl. 
\end{equation}
We look at three specific cases:

\noindent (i) (Isolated dislocation)
For an isolated dislocation at $O$ with Burgers vector $\boldsymbol{b}\in \mathbb{R}^2$, $\boldsymbol{A}=\boldsymbol{b}\delta_O$. The representation \eqref{represdrivforce} has only one non-trivial term such that
$\boldsymbol{F}^I=(\boldsymbol{\sigma}_2^o \boldsymbol{b})\times \boldsymbol{e}_3 \delta_O$, where $\boldsymbol{\sigma}_2^o = \boldsymbol{\sigma}_2(O)$. Hence,
\begin{equation}
\int_{\partial B_\epsilon} \boldsymbol{J}^I \boldsymbol{n} \dl=(\boldsymbol{\sigma}_2^o \boldsymbol{b})\times \boldsymbol{e}_3,
\end{equation}
which can be identified as the well known Peach-Koehler force acting on a dislocation in the presence of an external field $\boldsymbol{\sigma}_2$ \cite{eshelby1956continuum}. 

\noindent (ii) (Dislocation dipole) Let $\boldsymbol{A}=(\boldsymbol{b}\otimes \boldsymbol{v})\nabla \delta_O=\boldsymbol{b} \langle \nabla \delta_O, \boldsymbol{v}\rangle$, where $\boldsymbol{b}\in \mathbb{R}^2$ and $\boldsymbol{v}\in \mathbb{R}^2$ are constants. Such a dislocation density represents an isolated dislocation dipole at $O$. Indeed, consider two dislocations of equal and opposite Burgers vector $\boldsymbol{b}/h$, $h \in \mathbb{R}$, placed at points $O$ and $O+h\boldsymbol{v}$. Then, the dislocation density $\boldsymbol{A}_h=({\boldsymbol{b}}/{h})\delta_O-({\boldsymbol{b}}/{h})\delta_{O+h\boldsymbol{v}}$ takes the limiting value $\boldsymbol{A}$ as $h \to 0$. Using \eqref{GeneralizedDrivingForceonDistribution}, we obtain
\begin{equation}
\boldsymbol{F}^I=-\left(\nabla \boldsymbol{\sigma}_2^o(\boldsymbol{b}\otimes \boldsymbol{v}) \times \boldsymbol{e}_3 \right) \delta_O+ \left( \boldsymbol{\sigma}_2^o(\boldsymbol{b}\otimes \boldsymbol{v}) \nabla \delta_O \right) \times \boldsymbol{e}_3,
\end{equation}
where $\nabla \boldsymbol{\sigma}_2^o = \nabla \boldsymbol{\sigma}_2(O)$. The first term in the expression for $\boldsymbol{F}^I$ yields
\begin{equation}
\int_{\partial B_\epsilon} \boldsymbol{J}^I \boldsymbol{n} \dl=-\nabla \boldsymbol{\sigma}_2^o(\boldsymbol{b}\otimes \boldsymbol{v})\times \boldsymbol{e}_3 \label{forcedipole}
\end{equation}
as the force acting on the dislocation dipole at $O$ \cite{kroupa1966dislocation}. The second term in $\boldsymbol{F}^I$, on the other hand, represents a couple acting on the dislocation dipole. Such a couple can be interpreted as the generalized thermodynamic force which resists non-translational configurational changes in the dipole (e.g., the relative position of the two dislocations which constitute the dipole).  To elaborate, consider the dipole as a pair of dislocations with density $\boldsymbol{A}_h$. The corresponding generalized force is of the form $\boldsymbol{F}^I_h= \boldsymbol{f}_1 \delta_O+ \boldsymbol{f}_2 \delta_{O+h\boldsymbol{v}}$,
where $\boldsymbol{f}_1=({1}/{h}) ( (\boldsymbol{\sigma}_2^o \boldsymbol{b})\times \boldsymbol{e}_3)$ and $\boldsymbol{f}_2=-({1}/{h}) ( (\boldsymbol{\sigma}_2^h \boldsymbol{b})\times \boldsymbol{e}_3 )$, $\boldsymbol{\sigma}_2^h = \boldsymbol{\sigma}_2(O+h\boldsymbol{v})$. Then, $\boldsymbol{F}^I_h \to {\boldsymbol{F}^I}$ as $h\to 0$. In particular, $\boldsymbol{f}_1+\boldsymbol{f}_2$ converge to the net force on the dipole \eqref{forcedipole}.  The couple acting on the dipole is obtained as the limit of the net moment due to forces $\boldsymbol{f}_1$ and $\boldsymbol{f}_2$ on the individual dislocations. Note that when the external field $\boldsymbol{\sigma}_2$ is uniform, the force acting on the dipole vanishes but the acting couple is non-trivial. 

\noindent (iii) (Point defect)
Let $\boldsymbol{A}=({a}/{2})(\boldsymbol{e}_1 \otimes \boldsymbol{e}_2 - \boldsymbol{e}_2 \otimes \boldsymbol{e}_1)\nabla \delta_O$, where ${a}\in \mathbb{R}$ is constant. The dislocation density $\boldsymbol{A}$ can be interpreted either in terms of  two dislocation dipoles $\boldsymbol{A}_1=({a}/{2})(\boldsymbol{e}_1 \otimes \boldsymbol{e}_2)\nabla \delta_O$ and $\boldsymbol{A}_2=-({a}/{2})(\boldsymbol{e}_2\otimes \boldsymbol{e}_1)\nabla \delta_O$ or in terms of four dislocations (two from each of the dipole). The incompatibility associated with $\boldsymbol{A}$ is $N=-({a}/{2}) \Delta \delta_O$. The defect at $O$ can therefore be interpreted as a centre of dilation emerging from the presence of an isolated interstitial or vacancy. Using \eqref{GeneralizedDrivingForceonDistribution} we have
\begin{equation}
\boldsymbol{F}^I=a \nabla (\tr \boldsymbol{\sigma}_2^o)\delta_O+\frac{a}{2}\left( \boldsymbol{\sigma}_2^o(\boldsymbol{e}_1\otimes \boldsymbol{e}_2 - \boldsymbol{e}_2 \otimes \boldsymbol{e}_1)\nabla \delta_O\right)\times \boldsymbol{e}_3,
\end{equation}
where $\nabla (\tr \boldsymbol{\sigma}_2^o) = \nabla (\tr \boldsymbol{\sigma}_2)(O)$, whose first term yields
\begin{equation}
\int_{\partial B_\epsilon} \boldsymbol{J}^I \boldsymbol{n} \dl=a \nabla (\tr \boldsymbol{\sigma}_2^o)
\end{equation}
as the force acting on the centre of dilation  \cite{eshelby1956continuum}.  The second term in $\boldsymbol{F}^I$ is a couple acting on the defect in response to non-translational configurational changes associated with the defect.

\section{Concluding remarks} \label{conc}

\noindent \textit{Singular point at the end of a singular curve}: Most of our results in the preceding sections were obtained for fields (stresses, strains, forces, incompatibility, etc.) which were assumed to have singular support only at an isolated point $O \in \Omega$, while admitting smooth restrictions on $\Omega - O$. We can in fact extend our considerations to more general situations, where the singular support set includes additional points in $\Omega - O$, by treating the regular divergence and curl operators acting on the restricted fields as their distributional counterparts. Consider, for instance, $\boldsymbol{\sigma}\in \mathcal{D}'(\Omega,\Sym)$ and $\boldsymbol{B}\in \mathcal{D}'(\Omega,\mathbb{R}^2)$ such that $\singsupp(\boldsymbol{B})\subset \singsupp(\boldsymbol{\sigma})$ and $\singsupp(\boldsymbol{\sigma})$ is a set of $k$ discrete points in $\Omega$. The fields are smooth outside these $k$ points. Lemma~\ref{fullocalstressequib} can then be generalized by writing Equation~\eqref{restrictedBalance_1} for restrictions on $\Omega-\singsupp(\boldsymbol{\sigma})$ and Equations~\eqref{restrictedBalance_2}-\eqref{restrictedBalance_4} with respect to each of the $k$ points. As another example, we consider a stress field which, outside $O$, concentrates on a smooth curve ${S}$ with $O$ as one end point such that $\boldsymbol{\sigma}|_{\Omega-O}(\boldsymbol{\phi})=\int_{S} \left\langle \boldsymbol{\sigma}_1, \boldsymbol{\phi} \right\rangle \dl$, for all $\boldsymbol{\phi}\in \mathcal{D}(\Omega-O,\Lin)$, where $\boldsymbol{\sigma}_1:{S}\to \Sym$ is smooth and bounded. If $\deg(\boldsymbol{\sigma})<0$, then  $\boldsymbol{\sigma}$ is the unique extension of $\boldsymbol{\sigma}|_{\Omega-O}$, allowing us to write \cite{pandey2020topological}
 \begin{equation}
 \Div \boldsymbol{\sigma} (\boldsymbol{\psi})=\int_{S} \left\langle \div_S \boldsymbol{\sigma}_1+\kappa \boldsymbol{\sigma}_1 \boldsymbol{n}, \boldsymbol{\psi} \right\rangle \dl - \int_{S} \left\langle  \boldsymbol{\sigma}_1 \boldsymbol{n}  , \frac{\partial \boldsymbol{\psi}}{\partial \boldsymbol{n}} \right\rangle \dl-\boldsymbol{\sigma}_1 \boldsymbol{t}\delta_O,
 \end{equation}
for all $\boldsymbol{\psi} \in \mathcal{D}(\Omega,\mathbb{R}^2)$, where $\kappa$ is the curvature of the curve $S$, $\boldsymbol{n}$ is the unit normal to $S$, $\boldsymbol{t}$ is the unit tangent along $S$, $\div_S$ is the divergence along $S$, and ${\partial \boldsymbol{\psi}}/{\partial \boldsymbol{n}} = (\nabla \boldsymbol{\psi}) \boldsymbol{n}$. If we take $\boldsymbol{B}=\boldsymbol{B}_0\delta_O$, with constant $\boldsymbol{B}_0 \in \mathbb{R}^2$, then the local equilibrium equations can be obtained as $ \div_S \boldsymbol{\sigma}_1+\kappa \boldsymbol{\sigma}_1 \boldsymbol{n} = \boldsymbol{0}$ and  $\boldsymbol{\sigma}_1 \boldsymbol{n}= \boldsymbol{0}$ on $S-O$ and $\boldsymbol{\sigma}_1 \boldsymbol{t}=\boldsymbol{B}_0$ at $O$. As a third, and final, example, we consider an array of dislocations as the source of incompatibility. Accordingly, let $\boldsymbol{A}\in \mathcal{D}'(\Omega,\mathbb{R}^2)$ be such that $\singsupp(\boldsymbol{A})={S}$ and $\boldsymbol{A}|_{\Omega-O}(\boldsymbol{\phi})=\int_{S} \left\langle \boldsymbol{\alpha}_1, \boldsymbol{\phi} \right\rangle \dl$, for all $\boldsymbol{\phi}\in \mathcal{D}(\Omega-O,\mathbb{R}^2)$, where $\boldsymbol{\alpha}_1:{S}\to \mathbb{R}^2$ is smooth and bounded. If $\deg(\boldsymbol{A})<0$, then  $\boldsymbol{A}$ is the unique extension of $\boldsymbol{A}|_{\Omega-O}$. We can then obtain an expression for the incompatibility field $N = \Curl \boldsymbol{A}$ as
\begin{equation}
N ({\psi})=\int_{S} \left( -\kappa\langle \boldsymbol{\alpha}_1,\boldsymbol{t} \rangle+\left\langle (\nabla\boldsymbol{\alpha}_1)\boldsymbol{t},\boldsymbol{n} \right\rangle  \right)\psi \dl - \int_{S} \langle \boldsymbol{\alpha}_1,\boldsymbol{t} \rangle \frac{\partial {\psi}}{\partial \boldsymbol{n}}  \dl-\langle\boldsymbol{\alpha}_1 ,\boldsymbol{n}\rangle \delta_O,
\end{equation}
for all ${\psi} \in \mathcal{D}(\Omega)$, where ${\partial {\psi}}/{\partial \boldsymbol{n}} = \langle \nabla{\psi}, \boldsymbol{n} \rangle$. In particular if $S$ is a straight line and $\boldsymbol{\alpha}_1=\alpha_1 \boldsymbol{n}$, where $\alpha_1 \in \mathbb{R}$ is a constant,  then the incompatibility field is point supported with $N=-\alpha_1\delta_O$.

\noindent \textit{Multiply connected domains}:
To characterize the incompatibility of a strain field $\boldsymbol{E}\in \mathcal{D}'(\Omega,\Sym)$ over a multiply connected region $\Omega \subset \mathbb{R}^2$, the local condition $\Curl \Curl \boldsymbol{E}=N$ has to be supplemented by a topological (non-local) condition. For example, a smooth strain field $\boldsymbol{E}_0 \in C^{\infty}(\Omega-O,\Sym)$ over a doubly connected domain $\Omega-O$, which satisfies the local compatibility condition in $\Omega-O$, i.e., $\curl \curl \boldsymbol{E}_0=0$, can still be incompatible in the sense of the condition
 \begin{equation}
 \int_{L} \{ (\boldsymbol{E}_0)(\boldsymbol{y})   + (\boldsymbol{y}-\boldsymbol{x})\times \curl ((\boldsymbol{E}_0)(\boldsymbol{y})) \} \dy= -a_1\boldsymbol{e}_1+a_2\boldsymbol{e}_2 + a_3 \boldsymbol{e}_3\times \boldsymbol{x},
 \end{equation}
 where $L \subset \Omega-O$ is an arbitrary closed loop enclosing $O$ and $a_1,a_2,a_3 \in \mathbb{R}$ are sources of incompatibility located outside the domain at $O$ \cite[\S 156A]{Love}. Let $\boldsymbol{E}\in \mathcal{A}(\Omega,\Sym)$ be an arbitrary extension of $\boldsymbol{E}_0$ over the simply connected $\Omega$. The incompatibility field $N=\Curl \Curl \boldsymbol{E}$ is then point supported at $O$ such that $N^{(0,1)}=a_1,$ $N^{(1,0)}=a_2$, and $N^{(0,0)}=a_3$ (Lemma~\ref{IncompatibilityCeasroIntegral}). The incompatibility is therefore due to a dislocation and a disclination defect placed outside the body $\Omega-O$ at $O$. The components $N^{\alpha}$, for $|\alpha|>1$, are not fixed for an arbitrary extension of $\boldsymbol{E}_0$; they will differ for different extensions. Hence, it is meaningless to place higher-order defects, like dislocation dipole, extra matter, etc., outside the doubly connected domain $\Omega-O$. In fact, if we consider a higher order defect at $O$, such that  $N^{\alpha}=0$ for $|\alpha| \leq 1$, then the strain incompatibility condition (in terms of $\boldsymbol{E}$) does not  impose any necessary restrictions on $\boldsymbol{E}_0$ (Lemma~\ref{CurlCurlDipoleSource}).

\noindent \textit{Three-dimensional domains}: One of the main concerns of our work, which was restricted to planar domains, was to establish the sufficiency or insufficiency of the smooth maps, away from the point of singularity, in completely characterizing the equilibrium of singular stress fields, in the presence of singular body forces, and the compatibility/incompatibility of singular strain fields, in the presence of singular incompatible fields (arising from singular defect distributions). These results extend naturally to  fields with point singularities in a 3D domain. This is essentially due to the fact that the notion of the degree of distribution, with respect to the singular point $O$, and the uniqueness and existence of extensions of a distribution in $\Omega-O$ to a distribution in $\Omega$ are both given for $\mathbb{R}^n$ with arbitrary $n$. In a 3D domain $\Omega \subset \mathbb{R}^3$, however, fields which are singular on a curve, $L\subset \Omega$,  have to be necessarily considered to model isolated line defects such as dislocations and  disclinations \cite{scala2016currents}.  We will need results on the existence and uniqueness of extensions of a distribution on $\Omega-L$ to a distribution in $\Omega$ in order to formulate assumptions on the fields under which the smooth maps, away the curve of singularity, are sufficient to characterize the respective fields and the associated governing equations. Such results, which are currently unavailable, will provide the necessary foundations for developing a complete 3D theory incorporating point and line singularities.

\section*{Ackowledgement}
AG acknowledges the financial support from SERB (DST) Grant No. CRG/2018/002873 titled ``Micromechanics of Defects in Thin Elastic Structures". 

\appendix

\section{Proof of the existence of extension in Lemma~\ref{ExistenceUniquenessExtensionLemma}}
\label{appExtProof}

We first establish the existence result  when $\deg(T_0)<0$. Consider $\vartheta \in \mathcal{D}(\Omega)$ such that $\supp(\vartheta) \subset {B}_{r}$ and $\vartheta(\boldsymbol{x}) =1$ for all $\boldsymbol{x} \in {B}_{r/2}$. Given $\lambda>1$,  define $\vartheta_{\lambda} \in \mathcal{D}(\Omega)$ as $\vartheta_{\lambda} (\boldsymbol{x})=\vartheta(\lambda \boldsymbol{x})$. Hence, $\supp(\vartheta_\lambda) \subset B_{r/\lambda}$ and $\vartheta_{\lambda}(\boldsymbol{x}) =1$ for all $\boldsymbol{x} \in {B}_{{r}/{2\lambda}}$. For any $\phi \in \mathcal{D}(\Omega),$ $(1- \vartheta_{2^{j}})\phi \in \mathcal{D}(\Omega-O).$ We consider the sequence of distributions ${T}^j\in \mathcal{D}'(\Omega)$ as ${T}^j= (1- \vartheta_{2^{j}}){T_0}$. Hence, $T^j (\phi)={T_0}((1- \vartheta_{2^{j}})\phi)$
for all $\phi \in \mathcal{D}(\Omega).$ For any $\phi \in \mathcal{D}(\Omega),$ we have
$(T^{j+1}-T^j) (\phi)=(\phi T_0)(-\vartheta_{2^{j+1}}+ \vartheta_{2^{j}}) 
= 2^{-nj} \left((\phi T_0)|_{{B}_r} \right)_{2^{-j}} (\vartheta-\vartheta_2)$.
Given $\deg(T_0)<0$, for $k \in \mathbb{R}$ such that $\sd(T_0)< k <n$, we can conclude that
\begin{equation}
\lim_{j\to\infty} 2^{-nj} \left((\phi T_0)|_{\mathcal{B}_r} \right)_{2^{-j}} (\vartheta-\vartheta_2)=0.
\end{equation}
Also, there exists $j_0\in \mathbb{N}$ and $c_0\in \mathbb{R}$ such that, for any $j>j_0$, $(T^{j+1}-T^j) (\phi) < c_0 2^{j(k-n)}$.
Hence the sequence $T^j (\phi)$ is a Cauchy sequence and we can define a distribution $T\in \mathcal{D}'(\Omega)$ as
$T(\phi)=\lim_{j\to \infty} T^j (\phi)$.
Clearly, $T\in \mathcal{D}'(\Omega)$ such that $T|_{\Omega-O}=T_0$. It can be shown that $\deg(T)=\deg(T_0)$ \cite{brunetti2000microlocal}. Further, if $T_1 \in \mathcal{D}'(\Omega)$ is another extension of $T_0$ then $(T_1-T) \in \mathcal{E}(\Omega)$. Consequently $\deg(T_1)\geq 0$ if $(T_1-T)\neq 0$ (Lemma~\ref{RepresentationELemma}). This is a contradiction. Hence $T$ is the unique extension of $T_0$ such that $\deg(T)=\deg(T_0)$. Next, we consider the case when $\deg(T_0)  \geq 0$. Let $\rho$ be the greatest integer smaller (or equal) than $\deg(T_0)$ and let $x^{\alpha}={x_1}^{\alpha_1}{x_2}^{\alpha_2}\dots{x_n}^{\alpha_n}$ for any $\alpha=(\alpha_1,\alpha_2\dots\alpha_n).$ The function $\phi \in \mathcal{D}(\Omega)$ can be uniquely decomposed as
\begin{equation}
\phi=\sum_{|\alpha| \leq \rho} w^{\alpha} \partial^{\alpha} \phi (O)+\sum_{|\alpha|=\rho+1} x^{\alpha} \psi_{\alpha}
\end{equation}
where $\psi_{\alpha} \in \mathcal{D}(\Omega)$ and $w^{\alpha}\in \mathcal{D}(\Omega)$ is as introduced in Lemma~\ref{RepresentationELemma}.
 We have $\deg(x^{\alpha}T_0) <0$ for $|\alpha|=\rho+1$ \cite{brunetti2000microlocal}. Let $x^{\alpha}T \in \mathcal{D}'(\Omega)$ be the unique extension of $x^{\alpha}T_0 \in \mathcal{D}'(\Omega-O)$ such that $\deg(x^{\alpha}T)=\deg(x^{\alpha}T_0).$ The distribution $T \in \mathcal{D}'(\Omega)$, defined as
$T(\phi)=\sum_{|\alpha|=\rho+1} x^{\alpha}T(\psi_\alpha)$ satisfies $T|_{\Omega-O}=T_0$. It can be shown that $\deg(T)=\deg(T_0)$ \cite{brunetti2000microlocal}. Let $T_1 \in \mathcal{D}'(\Omega)$ be an extension of $T_0$ such that $\deg(T_1)=\deg(T_0)$. Then $(T_1-T) \in \mathcal{E}(\Omega)$ with $\deg(T_1-T)\leq\deg(T_0)$, and we can write $T_1-T=\Sigma_{|\alpha| \leq \deg({T_0})} q^{\alpha} \partial^\alpha \delta_O,$ where $q^{\alpha}\in \mathbb{R}$ (Lemma~\ref{RepresentationELemma}). 

%\begin{quote}
%\section{Introduction}

\medskip
 
\bibliographystyle{plain}
\bibliography{ref}

  \end{document}